\input graphicx



\def \bu { \hskip -.05in}
\def \LR {\leftrightarrow}
\def \ul {\underline}
\def\multi#1{\vbox{\baselineskip=0pt\halign{\hfil$\scriptstyle\vphantom{(_)}##$\hfil\cr#1\crcr}}}

\def \BN {{\bf N}}

\def \BF {{\bf F}}
\def \BZ {{\bf Z}}
\def \BF {{\bf F}}

\def\tttt #1{{\textstyle{#1} }}

\def \page{{\vfill\supereject}}

\def \magstep#1 {\ifcase#1 1000\or 1200\or 1440\or 1728\or 2074\or 2488\fi\relax}

 \font\Ch=msbm8

\def\Q{\hbox{\Ch Q}}

\def\la{{\lambda}}

\def \BF {{\bf F}}

\overfullrule=0pt
\baselineskip 14pt
\settabs 8 \columns
\+ \hfill&\hfill \hfill & \hfill \hfill   & \hfill \hfill  & \hfill \hfill   &
 & \hfill \hfill\cr
\parindent=.5truein
\hfuzz=1.44182pt
\hsize 6.5truein
\font\ita=cmssi10  
\font\small=cmr6
\font\title=cmbx10 scaled\magstep2
\font\normal=cmr10 
\font\small=cmr6

\font\bol=cmbx12

\def\sig{\sigma}

\def \-> {\rightarrow}
\def\LL{\big\langle}

\def\RR {\big\rangle}

\def\DD {\Delta}

\def\OM {\Omega}

\def\la {\lambda}
\def\La {\Lambda}
\def \RA {\rightarrow}

\def\xon {x_1,x_2,\ldots ,x_n}

\def \sas {\vskip .06truein}
\def\sa{{\vskip .125truein}}

\def\sap{{\vskip .25truein}}
\def\sapp {{\vskip .5truein}}

\def \eee {\epsilon}

\def\aa {\alpha}

\def \ses {\enskip = \enskip}
\def \sps {\, + \,}

\def \sms {\, - \,}

\def \scs {\, , \,}
\def \ess {\enskip}

\def \ssp {\hskip .25em}
\def \bigsp {\hskip .5truein}
\def \part {\vdash}

\def \DD {\Delta}

\normal

\vsize=9truein
\sap
\def\today{\ifcase\month\or
January\or February\or March\or April\or may\or June\or
July\or August\or September\or October\or November\or
December\fi
\space\number\day, \number\year}

\headline={   
\small  $\ess\ess\ess\ess\ess$ \hfill  
\hfill$\ess\ess\ess$\ess\ess\ess\ess\ess\ess$ \folio }
 \footline={\hfil}

\def \OM {{\Omega}}

\def \RA {{ \rightarrow }}

\def \BV {{\bf V}}

\def \BQ {{\Q}}

\def \II {1}

\def \TH {{\tilde H}}

\def\xon {x_1,x_2,\ldots ,x_n}

\def \BQ {{\Q}}

\def \TH {{\tilde H}}

\font\small=cmr6
\def \scs {\ssp , \ssp}
\def \ess {\enskip}
\def \ssp {\hskip .25em}
\def \bigsp {\hskip .5truein}
\def \part {\vdash}
\font\title=cmbx10 scaled\magstep2
\font\normal=cmr10 
\def\today{\ifcase\month\or
January\or February\or March\or April\or may\or June\or
July\or August\or September\or October\or November\or
December\fi
\space\number\day, \number\year}
\headline={\small  
\today    \small \hfill    \hfill$\ess\ess\ess\ess\ess\ess\ess\ess\ess$ \folio }
\null
\vskip -.5in

\def \II {{\rm I}}

\vskip15pt
\def \Be {{\bf e}}
\def \BF {{\bf F}}

\def \BC {{\bf C}}

\def \CPF {{\cal PF}_n}
\def \TS {\textstyle}

\def \and {\ess\ess\ess\hbox{and}\ess\ess}

\centerline{\bol A new Plethystic Symmetric Function Operator}
\centerline{\bf and }
\centerline{\bol The rational Compositional Shuffle Conjecture at t=1/q}
\centerline{\bf by}

\centerline{\bf A.M. Garsia, E. Leven, N. Wallach and G. Xin}
\sa

\noindent{\bf Abstract.}  Our main result here is that the specialization at $t=1/q$ of the $Q_{km,kn}$ operators studied in [4]   may be given  a very simple plethystic form. This discovery 
yields elementary  and direct derivations of several identities relating  these operators at $t=1/q$ to the Rational  
Compositional Shuffle conjecture of [3]. In particular we show that
if $m,n $  and  $k$  are positive integers and  $(m,n)$ is a coprime pair 
then 
$$
q^{(km-1)(kn-1)+k-1\over 2} Q_{km,kn}(-1)^{kn}\Big|_{t=1/q}\ses \tttt{[k]_q\over [km]_q} e_{km}\big[ X[km]_q\big]
\eqno {\rm A}.1
$$
where as customarily,  for any integer $s\ge 0$ and indeterminate $u$ we set $[s]_u=1+u+\cdots +u^{s-1}$.
We also show that the symmetric polynomial on the right hand side is always Schur positive. Moreover,  using the Rational Compositional Shuffle conjecture, we  derive a precise formula expressing this polynomial  in terms of Parking functions in  the $km\times kn$ lattice rectangle. 
\sa

\noindent{\bol Introduction}

The specializations at $t=1/q$ of all the Shuffle conjectures (including the classical cases) are still open to this date. What makes this specialization particularly fascinating is that both sides of
the stated identities have combinatorial interpretations. 
Nevertheless proving these identities is quite challenging even in the simplest cases. For instance  from the Rational Shuffle Conjecture we can easily derive the following identity,
for any coprime pair $(m,n)$.
$$
\sum_{D\in {\cal D}_{m,n}}q^{coarea(D) +dinv(D)}
\ses {1\over [m]_q }\left[ {m+n-1\atop n}\right]_q
\eqno \II.1
$$ 
where the sum is over Dyck paths in the $m\times n$ lattice rectangle, $coarea(D)$ gives the number of lattice squares 
above the path and $dinv(D)$ is a Dyck path statistic that can also be given a relatively simple geometric construction. The identity 
obtained by setting $q=1$ in I.1 is an immediate consequence of the  Cyclic Lemma, which suggests that this classical result
may have a natural $q$-analogue. The investigations that yielded the present results have been directed towards  the necessity of  
giving a concrete setting to a variety of identities stated or implied in recent work by  the Algebraic Geometers, most particularly
in the papers Burban-Schiffmann [5] and  Schiffmann-Vasserot [27]. Unfortunately most of this 
work appears in  language that  requires considerable algebraic Geometrical background. We have been privileged to have had the statements of
some of these results translated in a language that we could understand  by Eugene Gorsky and Andrey Negut. Many of the theorems we prove here have their  origin in this algebraic geometrical literature. Our contribution is to provide proofs
that are accessible to the algebraic combinatorial audience.
We hope that in doing so,   the new results we obtain may be conducive to progress in this most challenging area  of Algebraic
Combinatorics.

We will be dealing here with an algebra $\cal A$ of linear operators acting on the space $\La $ of symmetric function in an infinite alphabet $X=\{x_1,x_2,x_3,\ldots\}$ with coefficients in the field $\BQ(q,t)$ of rational functions in the two indeterminates $q,t$. Given a symmetric function $F[X]\in \La$, it will be convenient to denote by ``$\ul F$'' the operator ``{\ita multiplication by $F[X]$}''. As customary we will denote by ``$F^\perp$'' the operator dual of $\ul F$ with respect to the classical Hall scalar product of symmetric functions.

For a coprime pair $(m,n)$ the  $Q_{m,n}$ operators have an elementary definition, which as far as we understand, is due to 
Burban-Schiffmann in [5].  By taking a lattice point $(a,b)$  in the $m\times n$ rectangle that is closest and below the segment
$(0,0)\LR (m,n)$ and setting $(c,d)=(m,n)-(a,b)$
we obtain a decomposition
$$
(m,n)=(a,b)+(c,d) 
\eqno \II.2
$$
which here and after will be referred to as ``{\ita Split(m,n)}''.
Since the co-primality of the pair $(m,n)$ forces both pairs $(a,b)$ and $(b,c)$ to be coprime, we can recursively set
$$
Q_{m,n}\ses \tttt{1\over M}[Q_{c,d},Q_{a,b}]\scs
\bigsp (\hbox{$M=(1-t)(1-q)$})
\eqno \II.3
$$
with base cases 
$$
Q_{0,1}=  -\ul e_1
\ess\ess\ess\hbox{and }\ess\ess\ess
Q_{1,0}= D_0
\eqno \II.4
$$
where, $e_1$ is the customary elementary symmetric function and $D_0$ belongs to a family of operators $\{D_k\}_{k\in \BZ}$,
introduced  in [12], and defined by setting for $F[X]\in \La$
$$
D_kF[X]\ses F[X+\tttt{M\over z}]\sum_{r\ge 0}(-z)^re_r[X]
\Big|_{z^k}
\eqno \II.5
$$
In dealing with the present subject it is absolutely indispensable to 
make use of plethystic notation. Readers that are not familiar with this  device are referred to [14] for an introduction to its use. The readers will  also find in [3] and [4],  elementary proofs of all the auxiliary results that we will need in this writing. Those  papers were  written precisely to render this subject accessible to the algebraic combinatorial audience in a completely self contained manner.
In particular it is shown in [4] that to compute the action of an operator $Q_{m,n}$ we do not need to recurse to the base cases in I.4,  but rather use as a short cut  the identities
$$
Q_{1,k}\ses D_k
\eqno \II.6
$$
It should be mentioned that the original identities justifying the use of
this short cut were first given  in [12].

The definition of the operator $Q_{u,v}$ for a non coprime pair $ (u,v)$ relies on a truly amazing
property of the algebra generated by the operators $D_k$, certainly noticed in [5] and possibly in other algebraic geometrical   
literature. To state it it will be convenient to 
write a non-coprime pair in the form 
$(u,v)=(km,kn)$ with $(m,n)$  coprime 
and $k> 1$. This given, we can recursively define the operator $Q_{km,kn}$
by choosing any of the lattice points $(a,b)$ in the rectangle $km\times kn$
that are strictly below and closest to the segment $(0,0)\LR (km,kn)$
then set  
$$
Q_{km,kn}\ses \tttt{1\over M} [Q_{km-a,kn-b},Q_{a,b}].
\eqno \II.7
$$
This  definition is made possible  because, the choice of $(a,b)$   
forces both $(a,b)$ and $(km-a,kn-b)$ to be  coprime. Moreover,
all the  operators resulting from such a choice  of $(a,b)$ can be shown to
act identically on symmetric functions. 

Another fundamental fact discovered by the Algebraic Geometers is that the $Q$ operators indexed by collinear vectors do commute. More precisely for any coprime pair $m,n$ and any two integers $k,h$ we have
$$
[Q_{km,kn},Q_{hm,hn}]\ses 0.
\eqno \II.8
$$
Elementary, but by no means simple,  proofs of all these properties  
are given in [4]. The complexity of these proofs  is   due to the recursive nature of the definition in I.3.  Our discovery here is that with the specialization 
$t=1/q$ we can obtain several explicit identities from which these fundamental properties are immediate.
\sas

More precisely let $D_{u,v}$ denote  the   operator whose action on  the symmetric function $F[X]\in \La$ is defined by setting
$$
D_{u,v}F[X]\ses F\big[X+ \tttt{M[u]_q/z}\big]\sum_{r\ge 0}
(-z)^r e_r\big[ [u]_tX\big]\Big|_{z^v},
\eqno \II.9
$$
where  we  must set here $t=1/q$.
\sas

\noindent{\bol Theorem I.1}

{\ita If  $a,b,c,d,u,v$ are any integers related by the  vector identity
$
(a,b)+(c,d)= (u,v)
$
we have for non vanishing $a,c,u$}
$$
\tttt{1\over M}\big[D_{c,d},D_{a,b } \big]\Big|_{t=1/q}\ses q^{1+bc} {[a]_q[c]_q\over [u]_q}
{1-q^{da-bc}\over 1-q}D_{u,v}\Big|_{t=1/q}
\eqno \II.11
$$

This identity has the following  immediate corollary 
\sas

\noindent{\bol Theorem I.2}  

{\ita For any coprime pair $(m,n)$ and $k\ge 1$
we have}
$$
q^{(km-1)(kn-1)/2+( k-1)/2}  Q_{km,kn} \Big|_{t=1/q}= 
q^{(km-1)kn} 
\tttt{[k]_q\over [km]_q}      D _{km,kn}. 
\eqno \II.12
$$

The two identities in I.11 and I.12 have variety of   consequences. For instance we can immediately see from I.11 that the collinearity of $(a,b)$  and  $(c,d)$  implies that $ D_{a,b}$   and
$D_{c,d}$ commute. We thus obtain a much simpler proof of
this commutativity result for the $Q_{u,v}$ operators when $t=1/q$. Another immediate consequence of I.11 is that the algebra generated by the $D_k$ operators at $t=1/q$ is spanned by
the {\ita convex }  monomials in the $D_{u,v}$ operators. Here a monomial
$$
D_{u_1,v_1}D_{u_2,v_2}D_{u_3,v_3}\cdots D_{u_\l,v_\l},
$$
is called  {\ita convex }  if and only if we have
$$
\tttt {
{v_1\over u_1}\ge{v_2\over u_2}\ge \cdots \ge {v _\l\over u_\l}
}.
$$

To state an important consequence of I.12, we need some background. Let us recall that the classical Shuffle conjecture of 
Haglund et Al in [20] may be stated as the identity
$$
Q_{n+1,n}(-1)^n\ses \sum_{PF\in \CPF}t^{area(PF)}q^{dinv(PF)}
s_{pides(PF)}[X]
\eqno \II.13
$$
where the sum is over Parking Functions in the $n\times n$ lattice square,  $area(PF)$ and $dinv(PF)$ are Parking function statistics we will define later in a much more general context,
and $s_{pides(PF)}[X]$ denotes the Schur function indexed by the 
composition which gives the descent set of a permutation 
naturally associated to a Parking Function. In a recent paper
E. Gorsky and A. Negut formulated an infinite variety Shuffle conjectures, one for each coprime pair $(m,n)$.
They may be stated in a form similar to I.13, namely
$$
Q_{m,n}(-1)^n\ses \sum_{PF\in {\cal PF}_{m,n}}t^{area(PF)}q^{dinv(PF)}
s_{pides(PF)}[X]
\eqno \II.14
$$
where the sum is over Parking functions in the $m\times n$ 
lattice rectangle, and the parking function statistics occurring in
I.14 are highly non trivial modifications of the statistics involved in I.12.
Now, Theorem  I.2 has the following immediate corollary
 \sas
 
 \vbox{
 \noindent{\bol Theorem I.3}
 
{\ita For any coprime pair $(m,n)$ and $k\ge 1$ we have}
$$
q^{(km-1)(kn-1)/2+( k-1)/2}  Q_{km,kn}(-1)^{kn} \Big|_{t=1/q}\ses 
\tttt{[k]_q\over [km]_q}     e_{kn}\big[ X[km]_q\big]. 
\eqno \II.15
$$}

In particular, by combining I.15 with the Gorsky-Negut  conjectures at $t=1/q$  we obtain  the identity
$$
\tttt{1\over [m]_q}     e_{n}\big[ X[m]_q\big]
\ses
 \sum_{PF\in {\cal PF}_{m,n}} q^{coarea(PF)+dinv(PF)}
s_{pides(PF)}[X]
\eqno \II.16
$$
with   $coarea(PF)=(m-1)(n-1)/2-area(PF)$.

This given, we may ask if   the right hand side of I.15, can also be given a Parking Function interpretation. It turns out that this indeed the case.  More precisely we will show the following  
\sas

\noindent{\bol Theorem I.4}

{\ita Upon the validity of the extended Compositional Shuffle Conjecture in [3] it follows that
$$
\tttt{[k]_q\over [km]_q}     e_{kn}\big[ X[km]_q\big]
\ses
\sum_{PF\in {\cal PF}_{km,kn}} q^{coarea(PF)+dinv(PF)}
\big[ret(PF)\big]_q
s_{pides(PF)}[X]
\eqno \II.17
$$
where $ret(PF)$ is a statistic which indicates the height of the first 
return to the diagonal by the Dyck path of $PF$ in the $km\times kn$   lattice rectangle 
.}

The precise definitions of all the Parking Function statistics occurring in I.17 will be given in the sequel.

We must mention that it would follow from I.17 combined with the theory of LLT polynomials that the 
left hand side is a Schur positive symmetric polynomial. However,
we will show that 
this particular result can  be given a much more  elementary proof.
\sas

It is important to notice that   operators $D_{u,v}$ can be used for any integral values of $u$ and $v$. Now it follows from
I.12 for $m=1$ and $n=0$ that
$$
Q_{k,0}\Big|_{t=1/q}= D_{k,0}
\eqno \II.18
$$
It was known to the Algebraic Geometers that the family of operators $\{Q_{k,0}\}_{k\ge 1}$ have as a complete set of eigenfunctions the modified Macdonald basis 
$\{\TH_\mu[X;q,t]\}_\mu$ introduced in [11].
More precisely, we have
$$
Q_{k,0}\TH_\mu[X;q,t)\ses \Big(1-(1-t^k)(1-q^k)B_\mu(q^k,t^k)\Big )\TH_\mu[X;q,t],
\eqno \II.19
$$
where for a partition $\mu=(\mu_1,\mu_2,\ldots ,\mu_\l)$ we set
 $
B_\mu(q,t)= \sum_{i=1}^\l t^{i-1}\sum_{j=1}^{\mu_i}q^{j-1}
$.
Since it can be shown that the polynomial 
$\TH_\mu[X;q,t)$ specializes, at $t=1/q$  to a scalar multiple of $s_\mu\big[ \tttt{X\over 1-q}]$, it follows from I.19 and I.18
that we must also have
\sas

\noindent{\bol Theorem I.5}
$$
D_{k,0}s_\mu\big[ \tttt{X\over 1-q}]\ses
\Big(1-(1-q^{-k})(1-q^k)B_\mu(q^k,q^{-k})\Big )s_\mu\big[ \tttt{X\over 1-q}]
\eqno  \II.20
$$
Proving this identity  directly from the 
definition in I.9 leads to some highly non trivial combinatorial 
problems. However, with some effort, as we shall see, 
a  less direct but still entirely elementary path   to I.20 can actually be found. In fact this particular effort led to the discovery 
of the  following  formula  for the action of $D_{u,v}$
on the basis $\big\{s_\mu\big[ \tttt{X\over 1-q}]\big]_\mu$.
\sas

\noindent{\bol Theorem I.6}

{\ita For any $u,v>0$ and any partition $\mu$ we have
$$
 D_{u,v}s_\mu\big[\tttt{X\over 1-q}\big ]\ses
 (q^u-1)\sum_{i=1}^{  |\mu|+v} q^{up(\mu)_i+v -ui }
s_{p(\mu)+ve_i }\big[\tttt{X\over 1-q}\big ]
\eqno \II.21
$$
where $p(\mu)$ is the weak composition of length $ |\mu|+v $ obtained by adjoining zeros to the parts of $\mu$ and $e_i$ is the 
$i^{th }$ coordinate vector of length  $ |\mu|+v $. }
\sas

\noindent{\bol Remark I.1}

We should mention that there is another interesting by-product 
of our introduction of the operator $D_{u,v}$.
We learned from Eugene Gorsky  (see also section 6.10 of [16]) that in [29] it  is shown that for a suitable constant factor  $c_{m,n}(q)$
we have, for $(m,n)$ a coprime pair 
$$
Q_{m,n}\Big|_{t=1/q}\bu =\ssp  c_{m,n}(q) \nabla^{m\over n}\ul p_n  \nabla^{-{m\over n}}\Big|_{t=1/q} 
\eqno \II.22
$$
with  $\nabla$ the operator introduced in [2]. Now it turns out that one can easily derive I.22 from I.21 for $(u,v)=(m,n)$, directly from the definition of $\nabla$ given in [2]. 
\sas

This paper is divided into five sections. In the first section we 
give an elementary proof  of Theorem I.1. This type of proof has  been successfully used in various similar situations where
we needed a  straightforward  proof of an identity that was discovered by another path. In this section we also give proofs of Theorems  I.2 and I.3.

In section 2 we give the elementary argument  that proves the polynomiality and Schur positivity of the symmetric function in  I.15. 

In section 3 we give our Parking Function setting for the  symmetric function in  I.15. We also give there a simplified version of the Parking Function statistics that occur in the formulation of the Rational Compositional Shuffle Conjecture 
that take account of the most recent developments in this subject.

In section 4  we prove Theorems I.5 and I.6 and explore some of their consequences. What is interesting is that the  path to these proofs uses an argument that may be conducive to the discovery  of a variety other identities of  similar type.
 
In section 5 we give an outline of the theoretical steps  that led to the discovery of the operators $D_{u,v}$ and Theorem I.1. \sas

\noindent {\bf Acknowledgement}

We cannot  overemphasize the importance of the contributions to the very existence of this article by Eugene Gorsky and Andrey Negut.  The situation we find ourselves in  is that a variety of symmetric function constructs and combinatorial conjectures  that have been discovered and extensively studied by the Algebraic Combinatorists for more than 2 decades has been recently enriched in a truly spectacular fashion by discoveries of  the Algebraic Geometers, (see for instance
[5], [16], [17], [18], [9], [24], [26], [27], [28], [29]).
Unfortunately most of this information is written in a language that is  accessible only to people with extensive Algebraic Geometrical background. Over a period of nearly two years 
Gorsky and Negut provided us with extensive information  
about these developments in a language that we could understand. The present article is now the fourth in a series of articles ((3], [4], [9])  where our principal  aim has been to give  these 
contributions of the Algebraic Geometers  a  self contained
Algebraic Combinatorial setting, in the hope that this enrichment  of our subject will be conducive to further fruitful Algebraic Combinatorial developments.
\sa

\page

\noindent{\bol 1. Commutator properties of our new operators}
\sas

Our main  goal in this section is an elementary  proof of Theorem I.1 and its immediate Corollaries.
In a later section we will try to give a glimpse of the machinery that 
led to the discovery of the operators $D_{u,v}$ and yielded the original proof 
of this identity.
To give the reader an idea of the basic difference between these  
 two approaches,  we need to  recall a device which was extensively used in all previous work in the theory of Macdonald polynomials.
 We will refer to it as the ``$\OM$'' notation. The point of departure of plethystic substitutions is the operation of evaluating the  power symmetric function $p_k$
 at a formal power series $E=E(t_1,t_2,t_3,\ldots )$ containing an unlimited  number of indeterminates. We simply set
 $$
p_k[E]\ses =E(t_1^k,t_2^k,t_3^k,\ldots )
\eqno 1.1
$$ 
Since every symmetric function  can be expressed as a polynomial in the 
power functions, this definition allows to evaluate $F[E]$  for any given 
symmetric function $F$. This is what we  refer to as the plethystic substitution  of $E$ in $F$. In this vein we set
$$
\OM[E]\ses exp\bigg(\sum_{k\ge 1 }{p_k[E]\over k}\bigg).
\eqno 1.2
$$
Clearly, this definition implies that for any two expressions $A$ and $B$ we have
$$
\OM[A+B]=\OM[A ]\times \OM[ B]
\ess\ess\ess\ess\ess\ess\hbox{and}\ess\ess\ess\ess\ess\ess 
\OM[A-B]=\OM[A ]/ \OM[ B]
\eqno 1.3
$$
In particular it also follows from 1.1 and 1.2, that if 
$X=x_1+x_2+x_3+\cdots$ then
$$
\sum_{r\ge 0}(-z)^r e_r[X]\ses \OM[-zX].
\eqno 1.4
$$
Using this device the definition of the operators $D_k$ in I.5 can be rewritten in the form
$$
D_kF[X]\ses F[X+\tttt{M\over z}] \OM[-zX]
\Big|_{z^k}
\eqno 1.5
$$
Successive applications of two operators $D_a$ and $D_b$ to a symmetric function $F[X]$, in this notation, leads to the identities
$$
\eqalign{
D_bD_aF[X]
&\ses D_bF[X+\tttt{M\over z_1}] \OM[-z_1X]
\Big|_{z_1^a}
\cr
&\ses  F[X+\tttt{M\over z_1}+\tttt{M\over z_2}] \OM[-z_1(X+\tttt{M\over z_2})]
\OM[-z_2X]
\Big|_{z_1^az_2^b}
\cr
&\ses  F[X+\tttt{M\over z_1}+\tttt{M\over z_2}] 
\OM[-(z_1+z_2) X ]
\OM[-\tttt{Mz_1\over z_2})]
\Big|_{z_1^az_2^b}
\cr
&\ses  F[X+\tttt{M\over z_1}+\tttt{M\over z_2}] 
\OM[-(z_1+z_2) X ]
\tttt{ (1-  { z_1/z_2}) (1-  qt{ z_1 / z_2})
\over
(1-  t{ z_1/z_2}) (1-  q{ z_1/ z_2})}\tttt{1\over z_1^az_2^b}\Big|_{z_1^0z_2^0}
}
\eqno 1.6
$$
where the last equality results from multiple applications of the identities in 1.2.
By contrast if we carry out this calculation,
the way a computer would do it,   we would end up with the following
sequence of identities.
$$
\eqalign{
D_bD_aF[X]
&\ses D_b
\sum_{r_1\ge 0}F^{(r_1)}[X] \tttt{1\over z_1^{r_1}}
 \OM[-z_1X]
\Big|_{z_1^a}
\cr
&\ses D_b
\sum_{r_1\ge 0}F^{(r_1)}[X] 
(-1)^{r_1+a}e_{r_1+a}[X]
\cr
&\ses  
\sum_{r_1\ge 0}F^{(r_1)}[X+\tttt{M\over z_2}]  
(-1)^{r_1+a}e_{r_1+a}[X+\tttt{M\over z_2}]
 \OM[-z_2X]
\Big|_{ z_2^b}
\cr
&\ses  
\sum_{r_1,r_2\ge 0}
F^{(r_1,r_2)}[X] 
\tttt{1\over z_2^{r_2}} 
(-1)^{r_1+a}
\sum_{s=0}^{r_1+a}e_{r_1+a-s}[X]\tttt{1\over z_2^s}e_s[M] 
 \OM[-z_2X]
\Big|_{ z_2^b}
\cr
&\ses  
\sum_{r_1,r_2\ge 0}
F^{(r_1,r_2)}[X] 
(-1)^{r_1+a}
\sum_{s=0}^{r_1+a}e_{r_1+a-s}[X] e_s[M] 
(-1)^{r_2+s+b}e_{r_2+s+b}[X]
\cr}
\eqno 1.7
$$
where for convenience we have set
$$
F^{(r_1)}[X]=F[X+ M u_1]\Big|_{u_1^{r_1}}
\ess\ess\ess\ess\ess\ess\hbox{and}\ess\ess\ess\ess\ess\ess
F^{(r_1,r_2)}[X]=F[X+ M u_1+ M u_2]\Big|_{u_1^{r_1}u_2^{r_2}}
$$ 

We can see from this example that the second calculation of the action of the operator $D_bD_a$ is completely elementary and straight forward. But, in more complex situations, this approach  is not conducive to discovery but only to delivering 
the verification of an identity discovered by other means.  On the other hand 
the calculation of this action carried out in 1.6, in several significant instances, 
has led  to discovery and proof of surprising identities.
 Nevertheless,  we must add that due care must be taken in expressing the rational function, argument of the constant term
in 1.6, as an appropriate Laurent series in $z_1,z_2$. A systematic way of carrying this out in greater generality has been developed in [30]  and [31].

The proof of Theorem I.1 in this section will use the  approach
illustrated in 1.7.
The proof that follows the   approach in 1.6, led to the discover of the
operators $D_{u,v}$ and their commutator identities. This second proof  will be given in section 5.  
\sas

Recalling that the action of the operator $D_{u,v}$ is defined by setting
$$
D_{m,n}F[X]\ses F[X+[m]_q\tttt{M\over z}]\OM\big[ -z[m]_t X\big]\Big|_{z^n}
\eqno 1.8
$$
with the convention that $t=1/q$, we have

\noindent{\bol Theorem 1.1}

{\ita If  $a,b,c,d,m,n$ are any integers related by the  vector identity
$
(a,b)+(c,d)= (m,n)
$
we have for non vanishing $a,c,m$}
$$
\tttt{1\over M}\big[D_{c,d},D_{a,b } \big]\Big|_{t=1/q}\ses q^{1+bc} {[a]_q[c]_q\over [m]_q}
{1-q^{da-bc}\over 1-q}D_{m,n}\Big|_{t=1/q}.
\eqno 1.9
$$
\noindent{\bol Proof}

Using the notation   and the sequence of steps outlined in 1.7
we get
$$
\eqalign{
D_{a,b}F[X]
&= F[X+\tttt{M[a]_q\over z_1}]\OM[-z_1[a]_tX]\Big|_{z_1^b}
=\sum_{r_1\ge 0}F^{(r_1)}[X]\tttt{1\over z_1^{r_1}}\OM[-z_1[a]_tX]\Big|_{z_1^b}
\cr
&
=\sum_{r_1\ge 0}F^{(r_1)}[X](-1)^{r_1+b}e_{r_1+b}\big[[a]_tX \big]
\cr}
$$
and consequently
$$
\eqalign{
D_{c,d}D_{a,b}F[X]
&=\sum_{r_1,r_2\ge 0}F^{(r_1,r_2)}[X]
\tttt{(-1)^{r_1+b}\over z_2^{r_2}}
e_{r_1+b}\big[[a]_t(X+\tttt{M[c]_q\over z_2})  \big]\OM[-z_2[c]_tX]\Big|_{z_2^d}
\cr
&=\sum_{r_1,r_2\ge 0}F^{(r_1,r_2)}[X]
\tttt{(-1)^{r_1+b}\over z_2^{r_2}}
\sum_{s=0}^{r_1+b}
e_{r_1+b-s}\big[[a]_t X\big]
e_{s}\big[[a]_t M[c]_q  
 \big]
 \tttt{1\over z_2^s}
 \OM[-z_2[c]_tX]\Big|_{z_2^d}
 \ess\ess\ess\ess \ess\ess\ess\ess\ess   1.10
\cr
&=\sum_{r_1,r_2\ge 0}F^{(r_1,r_2)}[X]
\tttt{(-1)^{r_1+b} }
\sum_{s=0}^{r_1+b}
e_{r_1+b-s}\big[[a]_t X\big]
e_{s}\big[M [a]_t [c]_q  
 \big]
 (-1)^{r_2+d+s}e_{r_2+d+s}\big[ [ c]_tX]\big]
\cr}
$$
Now we can easy see that
$$
M [a]_t [c]_q\ses \tttt{(1-q)(1-1/q)(1-q^{-a})(1-q^c)\over (1-1/q)(1-q)}
\ses -q^{-a}(1-q^a)(1-q^c).
$$
Thus
$$
e_{s}\big[M [a]_t [c]_q   \big]
\ses \tttt{(-1)^s\over q^{as}}
h_s\big[ (1-q^a)(1-q^c) \big]
\ses \tttt{(-1)^s\over q^{as}}
(1-q^a)(1-q^c)\tttt{1-q^{sm}\over 1-q^m},
$$
and 1.10 becomes
$$
\eqalign{
\tttt{D_{c,d}D_{a,b}\over(1-q^a)(1-q^c)}F[X]
&=
\sum_{r_1,r_2\ge 0}F^{(r_1,r_2)}[X]
\tttt{(-1)^{r_1+b} }
\sum_{s=0}^{r_1+b}
e_{r_1+b-s}\big[[a]_t X\big]
\tttt{(-1)^s\over q^{as}}
\tttt{1-q^{sm}\over 1-q^m}
 (-1)^{r_2+d+s}e_{r_2+d+s}\big[ [ c]_tX]\big]
\cr
&=
\sum_{r_1,r_2\ge 0}F^{(r_1,r_2)}[X]
\tttt{(-1)^{r_1+r_2+n} }
\sum_{s=0}^{r_1+b}
e_{r_1+b-s}\big[[a]_t X\big]
\tttt{q^{-as}-q^{cs}\over 1-q^m}
 e_{r_2+d+s}\big[ [ c]_tX]\big]
\cr}
$$
Or better
$$
\tttt{(1-q^m)D_{c,d}D_{a,b}\over(1-q^a)(1-q^c)}F[X]
\ses A\sms B
\eqno 1.11
$$
With
$$
A\ses \sum_{r_1,r_2\ge 0}F^{(r_1,r_2)}[X]
\tttt{(-1)^{r_1+r_2+n} }
\sum_{s=0}^{r_1+b}
e_{r_1+b-s}\big[[a]_t X\big]
 q^{-as} 
 e_{r_2+d+s}\big[ [ c]_tX]\big]
$$
and
$$
B\ses \sum_{r_1,r_2\ge 0}F^{(r_1,r_2)}[X]
\tttt{(-1)^{r_1+r_2+n} }
\sum_{s=0}^{r_1+b}
e_{r_1+b-s}\big[[a]_t X\big]
 q^{cs} 
 e_{r_2+d+s}\big[ [ c]_tX]\big]
$$
Simple manipulation allow us to rewrite $A$ and $B$ in the more convenient forms
$$
\eqalign{
A
&\ses q^{ad} \sum_{r_1,r_2\ge 0}F^{(r_1,r_2)}[X]q^{a r_2 }
\tttt{(-1)^{r_1+r_2+n} }
\sum_{s=0}^{r_1+b}
e_{r_1+b-s}\big[[a]_t X\big] 
 e_{r_2+d+s}\big[q^{-a } [ c]_tX]\big]
\cr
}
$$
and
$$
B\ses q^{cb}\sum_{r_1,r_2\ge 0}F^{(r_1,r_2)}[X]q^{c r_1 }
\tttt{(-1)^{r_1+r_2+n} }
\sum_{s=0}^{r_1+b}
e_{r_1+b-s}\big[q^{-c}[a]_t X\big]
 e_{r_2+d+s}\big[ [ c]_tX]\big].
$$
Carrying out the interchanges  $a \leftrightarrow c $ and
$ b\leftrightarrow  d$ gives
$$
\eqalign{
\widetilde A
&\ses q^{cb} \sum_{r_1,r_2\ge 0}F^{(r_2,r_1)}[X]q^{c r_2 }
\tttt{(-1)^{r_1+r_2+n} }
\sum_{s=0}^{r_1+d}
e_{r_1+d-s}\big[[c]_t X\big] 
 e_{r_2+b+s}\big[q^{-c} [ a]_tX]\big]
\cr
}
$$
and
$$
\widetilde B\ses q^{ad}\sum_{r_1,r_2\ge 0}
F^{(r_2,r_1)}[X]q^{a r_1 }
\tttt{(-1)^{r_1+r_2+n} }
\sum_{s=0}^{r_1+d}
e_{r_1+d-s}\big[q^{-a}[c]_t X\big]
 e_{r_2+b+s}\big[ [ a]_tX]\big]
$$
Thus from 1.11 we derive that
$$
\tttt{(1-q^m)D_{a,b}D_{c,d}\over(1-q^a)(1-q^c)}F[X]
\ses  \widetilde A\sms \widetilde B
$$
and consequently
$$
\tttt{(1-q^m) [D_{c,d},D_{a,b} ]\over(1-q^a)(1-q^c)}F[X]
\ses   A\sms   B
\sms (\widetilde A\sms \widetilde B)
\ses (A+\widetilde B)\sms (B+\widetilde A)
\eqno 1.12
$$
Now we may rewrite $A$ by setting $u=r_1+b-s$ so
$0\le u\le r_1+b$ and $r_2+d+s=r_1+r_2+n-u$
obtaining
$$
\eqalign{
A
&\ses q^{ad} \sum_{r_1,r_2\ge 0}F^{(r_1,r_2)}[X]q^{a r_2 }
\tttt{(-1)^{r_1+r_2+n} }
\sum_{u=0}^{r_1+b}
e_u\big[[a]_t X\big] 
 e_{r_1+r_2+n-u}\big[q^{-a} [ c]_tX]\big]
\cr
}
$$
For $\widetilde B$ we set $u=r_2+b+s$ so
$r_2+b \le u\le r_1+r_2+n$ and  then make the switch $r_1\leftrightarrow r_2$
$$
\eqalign{
\widetilde B
&\ses q^{ad}\sum_{r_1,r_2\ge 0}
F^{(r_2,r_1)}[X]q^{a r_1 }
\tttt{(-1)^{r_1+r_2+n} }
\sum_{u=r_2+b}^{r_1+r_2+n}
 e_{u}\big[ [ a]_tX]\big]
 e_{r_1+r_2+n-u}\big[q^{-a}[c]_t X\big]
 \cr
 &\ses 
 q^{ad}\sum_{r_1,r_2\ge 0}
F^{(r_1,r_2)}[X]q^{a r_2 }
\tttt{(-1)^{r_1+r_2+n} }
\sum_{u=r_1+b}^{r_1+r_2+n}
 e_{u}\big[ [ a]_tX]\big]
 e_{r_1+r_2+n-u}\big[q^{-a}[c]_t X\big]
 \cr
 }
$$
This gives
$$
\eqalign{
A+\widetilde B
&\ses
 q^{ad}\sum_{r_1,r_2\ge 0}
F^{(r_1,r_2)}[X]q^{a r_2 }
\tttt{(-1)^{r_1+r_2+n} }
\sum_{u=0}^{r_1+r_2+n}
 e_{u}\big[ [ a]_tX]\big]
 e_{r_1+r_2+n-u}\big[q^{-a}[c]_t X\big]
 \cr
 &\bigsp\bigsp \sps 
  \sum_{r_1,r_2\ge 0}
F^{(r_1,r_2)}[X] 
\tttt{(-1)^{r_1+r_2+n} }
e_{r_1+b}\big[ [ a]_tX]\big]
 e_{ r_2+d}\big[ [c]_t X\big]
 \cr
 &\ses
 q^{ad}\sum_{r_1,r_2\ge 0}
F^{(r_1,r_2)}[X]q^{a r_2 }
\tttt{(-1)^{r_1+r_2+n} }
 e_{r_1+r_2+n}\big[[m]_t X\big]
 \cr
 &\bigsp\bigsp \sps 
  \sum_{r_1,r_2\ge 0}
F^{(r_1,r_2)}[X] 
\tttt{(-1)^{r_1+r_2+n} }
e_{r_1+b}\big[ [ a]_tX]\big]
 e_{ r_2+d}\big[ [c]_t X\big]
 \cr
  &\ses
 q^{ad}\sum_{r_1,r_2\ge 0}
F^{(r_1,r_2)}[X]\tttt{q^{a r_2 }\over z^{r_1+r_2}}
\OM\big[ -[m]_tX\big]\Big|_{z^n}
 \cr
 &\bigsp\bigsp \sps 
  \sum_{r_1,r_2\ge 0}
F^{(r_1,r_2)}[X] 
\tttt{(-1)^{r_1+r_2+n} }
e_{r_1+b}\big[ [ a]_tX]\big]
 e_{ r_2+d}\big[ [c]_t X\big]
\cr
  &\ses
 q^{ad}
F\big[X+\tttt{M[m]_q\over z} \big]
\OM\big[ -[m]_tX\big]\Big|_{z^n}
 \cr
 &\bigsp\bigsp \sps 
 \sum_{r_1,r_2\ge 0}
F^{(r_1,r_2)}[X] 
\tttt{(-1)^{r_1+r_2+n} }
e_{r_1+b}\big[ [ a]_tX]\big]
 e_{ r_2+d}\big[ [c]_t X\big].}
$$
With an entirely analogous  sequence of steps we obtain the identity
$$
\eqalign{
B+\widetilde A
&\ses
 q^{cb} 
 F\big[X+\tttt{M[m]_q\over z} \big]
 \OM\big[ -[m]_tX\big]\Big|_{z^n}
\cr
&
\bigsp\bigsp\sps
  \sum_{r_1,r_2\ge 0}F^{(r_1,r_2)}[X]  
  \tttt{(-1)^{r_1+r_2+n} }
e_{r_2+d}\big[[c]_t X\big] 
 e_{r_1 +b}\big[ [ a]_tX]\big].  
\cr}
$$
Thus   1.12  finally yields
$$
\tttt{(1-q^m) [D_{c,d},D_{a,b} ]\over(1-q^a)(1-q^c)}
[D_{c,d},D_{a,b} ]  F[X]
\ses 
(q^{ad}-q^{cb})F\big[X+\tttt{M[m]_q\over z} \big]
\OM\big[ -[m]_tX\big]\Big|_{z^n}.
$$
Proving the identity
$$
\tttt{1\over M} [D_{c,d},D_{a,b} ]F[X]
\ses 
\tttt{q[a]_q[c]_q
\over [m]_q   }
\tttt{q^{cb}-q^{ad}\over 1-q}F\big[X+\tttt{M[m]_q\over z} \big]
\OM\big[ -[m]_tX\big]\Big|_{z^n},
$$
which is just another way of writing 1.9.
\sa

In the remainder of this section we will derive a number of immediate consequences of the identity in 1.9. We will state and prove them as a succession of Corollaries, the  last two of which are simply restatements of  Theorems 
I.1 and I.2. For convenience we will here and after use the symbol 
``$Q^s_{u,v}$'' as  a short hand for   ``$Q _{u,v}\Big|_{t=1/q}$'' 
We will start with an auxiliary identity that we will use in various occasions.
\sas

\noindent{\bol Lemma 1.1}

{\ita  For integers $a,b,c,d$, we have}
  $$ 
  (a-1)(b+1)/2 +(c-1)(d+1)/2 +bc+1 = (a+c-1)(b+d+1)/2 +(bc-ad+1)/2.
  $$

\noindent{\bol Corollary 1.1}

{\ita If $m$ and $n$ are relatively prime, then we have}
 $$
 Q^s_{m,n} = q^{(n+1)(m-1)/2} D_{m,n}/[m]_q.
 $$
\noindent{\bol Proof}

We will proceed by induction on $m$.
The base case $m=1$ is easy:  we have
  $Q^s_{1,n}=D_n\Big|_{t=1/q}= D_{1,n}$, so the corollary holds in this case.

It is best to start by an example. Let us consider $Q^s_{2,n}={1\over M}[Q^s_{1,d},Q^s_{1,b}]$,
where $Split(2,n)=(1,b)+(1,d)$, with $d-b=1$. Then by 1.9
$$
\eqalign{
  Q^s_{2,n}&
  = {1\over M}[D_{1,d},D_{1,b}] = q^{b+1} \tttt{ [1]_q[1]_q\over [2]_q}  \tttt{(1-q^{d-b})\over 1-q}D_{2,n} \cr
 &= q^{(m-1)(n+1)/2} \tttt{1\over [2]} D_{2,n}.
\cr
}
$$
One more example. Let us say $Split(3,n)=(1,b)+(2,d)$ with $d-2b=1$. Then again by 1.9
$$
\eqalign{
  Q^s_{3,n}&= \tttt{1\over M} [Q^s_{2,d},Q^s_{1,b}] 
         = q^{(d+1)/2}  \tttt{1\over[2]} [D_{2,d},D_{1,b}] 
          =  q^{(d+1)/2} q^{2b+1}   \tttt{1\over[3]} D_{3,n}\cr
         &= q^{(n+1)} \tttt{1\over[3]} D_{3,n}. 
\cr
}
$$
 Assume the corollary holds for smaller $m$. Now suppose
  $$ Split(m,n)=(a,b)+(c,d). $$
That is we have  $a+c=m, b+d=n, ad-bc=1$. Then
$$
\eqalign{
  Q^s_{m,n} &=  \tttt{1\over M} [Q^s_{c,d}, Q^s_{a,b}] 
  \cr
&= q^{(c-1)(d+1)/2+(a-1)(b+1)/2} \tttt {1\over [a]_q [c]_q} 
  \tttt {1\over M} [D_{c,d},D_{a,b}] 
   \cr
&=q^{(c-1)(d+1)/2+(a-1)(b+1)/2} \cdot q^{bc+1}
 \tttt{1\over[a+c]_q}   
 \tttt {1-q^{ad-bc}\over 1-q}   D_{a+c,b+d}
 \cr
&= q^{(c-1)(d+1)/2+(a-1)(b+1)/2+bc+1} \tttt {1\over[m]_q}      D_{m,n}
\cr
&= q^{(m-1)(n+1)/2} \tttt {1\over[m]_q}D_{m,n}.
\cr
}
$$
This completes the induction and the proof.

\noindent{\bol Corollary  1.2}

{\ita
For any coprime pair $(m,n)$   we have}
$$
q^{(m-1) (n-1)/2}Q^s_{m,n} (-1)^n =   \tttt{1\over [m]_q}\,    e_n[ X [m]_q]. 
$$
\noindent{\bol Proof}

 $$
 \eqalign{
Q^s_{m,n} (-1)^n 
&= q^{(m-1) (n+1)/2}  \tttt {1\over [m]_q} D_{m,n} (-1)^n \cr 
&=q^{(m-1) (n+1)/2}   \tttt{1\over [m]_q}  (-1)^n \Omega{-z X [m]_t} z^{-n} \Big|_{z^0} \cr
&=q^{(m-1) (n+1)/2}   \tttt{1\over [m]_q}   e_n[ X [m]_t]  
\cr
&=q^{(m-1) (n+1)/2} q^{-(m-1)n}  \tttt {1\over [m]_q}   e_n[ X [m]_q]  
\cr
&=q^{-(m-1) (n-1)/2}   \tttt{1\over [m]_q}   e_n[ X [m]_q]. 
\cr
}
$$
\sas

\noindent{\bol Corollary  1.3}

{\ita For any coprime pair $(m,n)$ and $k\ge 1$   we have}
$$
Q^s_{km,kn} = q^{(km-1)(kn+1)/2-(k-1)/2}\tttt {[k]_q\over [km]_q} D_{km,kn}.
$$
\noindent{\bol Proof}

Suppose
 $ 
Split(m,n)=(a,b)+(c,d). 
 $
That is  $a+c=m, b+d=n, ad-bc=1$. It should be clear that $a$, $b$ are also relatively prime, 
we can choose $Split (km,kn)=(a,b)+((k-1)a+kc,(k-1)b+ kd),$ so that it is easily obtained by linear algebra that
$$ 
det \pmatrix{
         a & (k-1)a+kc \cr
         b & (k-1)b+kd \cr
       } =k ,
\quad \det\pmatrix{
         a+c & (k-1)a+kc \cr
         b+d & (k-1)b+kd \cr
       } =1 , 
$$
Thus $c'=(k-1) a+kc$ and $d'=(k-1)b+kd$  are relatively prime. We have
$$
\eqalign{
Q^s_{km,kn} 
&= \tttt {1\over M} [Q^s_{(k-1)a+kc,(k-1)b+kd}, Q^s_{a,b}] 
\cr
&= q^{(c'-1)(d'+1)/2+(a-1)(b+1)/2} \tttt {1\over [a]_q [c']_q}  \tttt{1\over M} [D_{c',d'},D_{a,b}]
\cr
&=q^{(c'-1)(d'+1)/2+(a-1)(b+1)/2} \cdot q^{bc'+1} \tttt{1\over [a+c']_q}   
\tttt{1-q^{ad'-bc'}\over 1-q}   D_{a+c',b+d'}
\cr
&= q^{(c'-1)(d'+1)/2+(a-1)(b+1)/2+bc'+1} \tttt{[k]_q\over[km]_q} D_{km,kn} \cr
&= q^{(km-1)(kn+1)/2-(k-1)/2} \tttt {[k]_q\over [km]_q}      D_{km,kn}.
\cr
}
$$
\sas

\noindent{\bol Corollary  1.4}

{\ita For any coprime pair $(m,n)$ and $k\ge 1$   we have}
  $$
  Q^s_{km,kn} (-1)^{kn} =q^{-(km-1) (kn-1)/2 -(k-1)/2 } \tttt  {[k]_q\over [km]_q}   e_n\big[ X [km]\big]. 
$$ 
\noindent{\bol Proof}

We have
$$
\eqalign{
Q^s_{km,kn} (-1)^{kn} 
&=  q^{(km-1)(kn+1)/2-(k-1)/2}  {[k]_q\over [km]_q}      D_{km,kn} (-1)^{kn}
\cr
&= q^{(km-1)(kn+1)/2-(k-1)/2}  {[k]_q\over [km]_q}      D_{km,kn} (-1)^{kn} (-1)^{kn} \Omega{-z X [km]_t}   \Big|_{z^{kn}}
\cr
&=q^{(km-1) (kn+1)/2-(k-1)/2}  {[k]_q\over  [km]_q}   e_n\big[ X [km]_t\big] 
 \cr
&=q^{(km-1) (kn+1)/2 -(k-1)/2 } q^{-(km-1)n}   {[k]_q\over [km]_q}   e_n\big[ X [km]_q\big] 
\cr
&=q^{-(km-1) (kn-1)/2 -(k-1)/2 }   {[k]_q\over  [km]_q}   e_n\big[ X [km]_q\big].
\cr
}
$$

\page

\noindent{\bol 2. Polynomiality and positivity.}

The main goal of this section is to prove that the quotient 
$$
{[k]_q\over [km]_q}     e_{kn}\big[ X[km]_q\big]
\eqno 2.1
$$
is a Schur positive symmetric polynomial.
This will be obtained by combining the next four auxiliary  propositions.

The first fact on which the proof is based is the following classical
result.
\sas

\noindent{\bol Proposition 2.1} 

{\ita For any  $n\ge 1$ we have to  matrices $\|c_{\la,\rho}\|_{\la,\rho\part n}$
and  $\|d_{\la,\rho}\|_{\la,\rho\part n}$ such that }
$$
a)\ess\ess 
s_\la= \sum_{\rho\part n}c_{\la,\rho}p_\rho
\ess\ess\ess\ess\ess\ess \hbox{and }\ess\ess\ess\ess\ess\ess
b)\ess\ess 
p_\rho= \sum_{\la\part n}d_{\la,\rho}s_\la
\eqno 2.2
$$
\noindent{\bol Proof}

Frobenius proved that  2.2 a) and b)  hold  with
$$
a)\ess\ess  c_{\la,\rho}= \chi^\la_\rho/z_\rho
\ess\ess\ess\ess\ess\ess \hbox{and }\ess\ess\ess\ess\ess\ess
d)\ess\ess  d_{\la,\rho}= \chi^\la_\rho 
\eqno 2.3
$$
with $\chi^\la_\rho$ the value of Young's  irreducible  $S_n$ character indexed by $\la$ at the conjugacy class indexed by $\rho$ 
and for $\rho=1^{\aa_1}2^{\aa_2}\ldots n^{\aa_n}$ we have
$$
z_\rho=1^{\aa_1}2^{\aa_2}\cdots n^{\aa_n}\aa_1!\aa_2!\cdots \aa_n!
$$

This given, our polynomiality result   can be stated as follows
\sas

\noindent{\bol Proposition 2.2}

{\ita If $(m,n)=1$ and $k\ge 1$ then for all $\la\part kn$ we have
$$
\tttt{[k]_q\over [km]_q}s_\la\big[[km]_q\big]\in  \BQ[q]
\eqno 2.4
$$ 
where for any integer $s\ge 0$ we set $[s]_q=1+q+\cdots +q^{s-1}$}

\noindent{\bol Proof}

Note that to show 2.4 we need only show that every root of 
$$
1+q+q^2+\cdots +q^{km-1}\ses 0
\eqno 2.5
$$
is a root of the polynomial
$$
 [k]_qs_\la\big[[km]_q\big].
$$
To show this we use 2.2 (for $n\RA kn$) and write for $\zeta$ a root of 2.5 
$$
[k]_\zeta s_\la\big[[km]_\zeta\big]= \sum_{\rho\part kn}
c_{\la,\rho} [k]_\zeta \ssp p_\rho \big[[km]_\zeta\big]
\eqno 2.6
$$
Thus if the left hand side of $2.6$ does not vanish we will necessarily have a $\rho\part kn$
such that
$$
[k]_\zeta \ssp p_\rho \big[[km]_\zeta\big]\neq 0
\eqno 2.7
$$
in particular if for some $r$ we have $\rho_r=i$ then   
$$
1+\zeta^i+\zeta^{2i}+\cdots \zeta^{ (km-1)i}\ses p_i\big[[km]_\zeta  \big]\neq 0
\eqno 2.8
$$
and since
$$
(1-q^i)\big(1+q^i+q^{2i}+\cdots q^{ (km-1)i}\big)\ses 1- q^{(km)i}
$$
from 2.8, and the fact that $\zeta$ is a root of 2.5,  we derive that 
$$
a)\ess\ess \zeta^i=1\ess\ess\ess\ess\ess\ess\&\ess\ess\ess\ess\ess\ess 
b)\ess\ess \zeta\neq 1
\eqno 2.9
$$
Now note that since the $\rho$ in 2.7 may be written in the form 
$\rho=\prod_{i=1}^{kn}i^{\aa_i}$ with $\sum_{i=1}^{kn}i \aa_i=kn$ 
then we must also have
$$
\zeta ^{kn}=\prod_{i=1}^{kn}(\zeta^{i})^{\aa_i}\ses 1
$$
Now the assumed primality of the pair $m,n$ gives that $k=gcd(km,kn) $ and this  
combined with the fact that $\zeta$ is a root of 2.5 forces
$$
\zeta^k=1
$$
but then b) of 2.9 yields 
$$
1+\zeta+\zeta^2+\cdots +\zeta^{k-1}=0
$$ 
which is in plain contradiction with 2.7. This contradiction forces every root of 2.5
to be a root of 
$$
 [k]_qs_\la\big[[km]_q\big]\ses 0 
 $$
as desired.
\sas

The next device we use  is the following well known fact
\sas

\noindent{\bol Proposition 2.3}

{\ita Any principal evaluation of a Schur Function is unimodal. More precisely for any  $\la\part n$ and any $m>1$ the polynomial
$$
s_\la\big[[m]_q\big]
\eqno 2.10
$$
is unimodal.}

\noindent{\bol Proof}

This   is exercise 4 page 137 of Macdonald book [25]. Since the solution in [25] is barely sketched, for sake of completeness we carry out Macdonald's exercise in full detail.
 Macdonald considers the evaluation
$$
s_\la\Big[
{x_1^{m}-x_2^{m}
\over
x_1-x_2
}
\Big]\ses s_\la\big[
s_{m-1}[x_1+x_2]
\big]
\eqno 2.11
$$
as a character of $GL_2[\BC]$. Using this he derives that for some weakly positive integer constants $c_{r_1,r_2}$ we must have the expansion
$$
 s_\la\big[
s_{m-1}[x_1+x_2]
\big]= \sum_{\multi{ r_1+r_2=d\cr r_1\ge r_2} }
c_{r_1,r_2}s_{[r_1,r_2]}[ x_1+x_2]
\eqno 2.12
$$
where for convenience we have set 
$$
d=(m-1)n
\eqno 2.13
$$
\noindent
Notice that we have
$$
\eqalign{
s_{[r_1,r_2]}[ x_1+x_2]
&\ses
x_1^{r_1}x_2^{r_2}+ x_1^{r_1-1}x_2^{r_2+1}+
\cdots +
 x_1^{r_2+1}x_2^{r_1-1}+x_1^{r_2}x_2^{r_1}
 \cr
&\ses
x_2^d
\big(
(\tttt{x_1\over x_2})^{r_1}+
(\tttt{x_1\over x_2})^{r_1-1}+
\cdots +
(\tttt{x_1\over x_2})^{r_2}
   \big)
  \ses 
x_2^d
\sum_{s=r_2}^{r_1}q^s
 \cr }
$$
where for convenience we have set $\tttt{x_1\over x_2}=q$.
 Thus with a slight change of notation we may rewrite 2.12 in the form
$$
 s_\la\big[
s_{m-1}[x_1+x_2]
\big]= x_2^d\sum_{ r_2=0}^{\lfloor d/2\rfloor  } 
c_{r_2}\sum_{s=0} ^{d}q^s
\chi\big(r_2\le s\le d-r_2\big)
=x_2^d\sum_{s=0}^d q^s\sum_{r_2=0}^{s \wedge (d-s)}c_{r_2}
\eqno 2.14
$$
from which the unimodality assertion immediately follows
by setting $x_2=1$ and $x_1=q$.
\sa

Our positivity result is a consequence of the following simple but powerful  fact
\sas

\noindent{\bol Proposition  2.4}

{\ita Let $g(q)=b_0+b_1q+\cdots + b_rq^r $   
and assume that, for $d=r+s$, the  polynomial 
$$
f(q)=(1+q+\cdots +q^s)g(q)=\sum_{\l=0}^dc_\l \, q^\l
\eqno 2.15
$$
is  unimodal  with peak at $p$
and non-negative coefficients. Then $g(q)$ also has non-negative coefficients. }

\noindent{\bol Proof}

We proceed by contradiction. Suppose some of the coefficients of $g(b)$ 
are negative. Since $c_0=c_d\ge0$ let $b_i$ and $b_j$ (with $0<i\le j<d$)    be the leftmost and rightmost negative coefficients of $g(q)$ respectively.
Now if $i\le p$ then
\vskip -.3in
$$
c_i=\sum_{u=0\vee(i-s)}^i b_u
\ess\ess\ess\ess\& \ess\ess\ess\ess
c_{i-1}=\sum_{u=0\vee(i-1-s)}^{i-1} b_u
$$
This gives
$$ 
c_i-c_{i-1}\ses b_i\sms \chi\big(i-1-s>0\big)\, b_{i-1-s}<0 
$$
a contradiction! 

If $r-j\le d-p$ then do the same argument for $\tilde f(q)=q^d f(1/q)$ and
$\tilde g(q)=q^rg(1/q)$. So we are left with $i>p$ and $r-j>d-p$. But that cannot happen since it implies that
 $\ess
i>s+j>j.
 $

As a corollary we obtain our desired goal.
\sas

\noindent{\bol Theorem 2.1}

{\ita For any coprime pair $(m,n)$ and any $k\ge 1$ we have that the symmetric function
$$
\tttt{[k]_q\over [km]_q}e_{kn}\big[X[km]_q\big]
\eqno 2.16
$$
is a Schur positive symmetric polynomial}

\noindent{\bol Proof}

The Cauchy formula gives
$$
\tttt{[k]_q\over [km]_q}e_{kn}\big[X[km]_q\big]\ses 
\sum_{\la\part kn}s_{\la '}[X]
\tttt{[k]_qs_\la \big[[km]_q\big]\over [km]_q }
$$
\noindent
Now we have proved that
 $
[k]_qs_\la \big[[km]_q\big]
 $
is divisible by $[km]_q$, we have also proved that 
$s_\la \big[[km]_q\big]$ is palindromic unimodal. Thus it follows from this that also $[k]_qs_\la \big[[km]_q\big]$
is palindromic unimodal, we can then apply Proposition 2.4 
with $f(q)=[k]_qs_\la \big[[km]_q\big]$ and 
$$
g(q)  \ses  {[k]_qs_\la \big[[km]_q\big]\over [km]_q }
$$
and conclude that $g(q)\in \BN[q]$, proving the 
 Schur positivity of the polynomial in 2.16.

\supereject

\noindent{\bol Remark 2.1}

The polynomiality of the symmetric function in 2.16 for the special case $k=1$ was first proved by Mark Haiman in [22]. The question
arose from the discovery in [10] that in the case $m=n+1$ 2.16 
(for $k=1$) is the Frobenius characteristic of an appropriate single grading of the  Diagonal Harmonic  Module of $S_n$. This also prompted Haiman to seek for some reason that justified this polynomial to be Schur positive. With remarkable foresight Haiman investigates the more general $(m,n)$ case and provided a mechanism for proving that the symmetric function
$$
e_n\big[ X[m]_q\big]\over
[m]_q
\eqno 2.17
$$
is a polynomial and Schur positive if and only if $(m,n)$ is a coprime pair. However, the Schur positivity was shown in [22] by constructing   
a quotient
$$
\BQ[\xon]/(e_1,f_2,\ldots ,f_n)
$$
with Frobenius characteristic the polynomial in 2.17, where 
 $e_1$ the ordinary elementary  symmetric function and
$
f_1,f_2,\ldots ,f_n
$
a sequence of polynomials satisfying the following properties 

\itemitem {(1)} {\ita each $f_i$ is homogeneous of degree $m$
}

\itemitem {(2)} {\ita they satisfy the identities}
 $\ess 
\sig f_i=  f_{\sig_i}
\ess {
(\hbox{\ita for  $1\le i\le n$ and all $\sig\in S_n$}
)
}
 $

\itemitem {(3)} $f_1+f_2+\cdots +f_n=0$

\itemitem {(4)} $e_1,f_2,\ldots. f_n$ are a regular sequence.

\noindent
A sequence $f_i$ satisfying (1),(2),(3),(4) was constructed by 
Hanspeter Kraft a few years later but never published.
More recently a very natural example of such a sequence  was discovered by Dunkl  in [6] and used later by Gorsky in  his work [16] on torus knots invariants. It follows from this and Mark Haiman result  the Schur positivity of the polynomial in 2.17 in the coprime case.

The challenge now is to construct an equally natural quotient with
Frobenius characteristic the polynomial in 2.16.
\sa  

\page

\noindent{\bol 3. A Parking Function  setting for our Frobenius Characteristics}
\sas

The main goal of this section is a proof of Theorem I.4. To carry this out we need to briefly review the statement of the Rational 
Compositional Shuffle Conjecture. We will start by  introducing  the symmetric function tools that are used in its formulation.

The  basic ingredient here is the identity
$$
Q_{0,k}\ses  
\tttt{qt\over qt-1}
\,  \ul h_{k}[X(1/qt-1)] 
\ess\ess\ess\ess
\hbox{(for all integers $k\ge 1$)}.
\eqno 3.1
$$
From this it follows that a basis for the subspace of operators 
of bi-degrees  $(0,n)$ for $n\ge 1$ is given by the collection
\vskip -.3in
$$
\big\{
Q_{0,\la}
\big\}_\la
\ses
\Big\{
\prod_{i=1}^{\l(\la)}Q_{0,\la_i}
\Big\}_\la
\eqno 3.2
$$
\vskip -.06in

\noindent
This fact, (see [3] for an elementary treatment), can be used 
 to construct an operator of bi-degree $(km,kn)$ for any 
 coprime pair $(m,n)$, any integer $k\ge 1$ 
 and any given homogeneous symmetric function $F[X]$ of degree
 $k$, by the following two steps.

\item {(1)} {\ita Compute the expansion}
 $$
 F\ses \sum_{\la\part k}c_\la \prod_{i=1}^{\l(\la)}Q_{ 0, \la_i},
\eqno 3.3
 $$

\vskip -.1in
 \item {(2)} {\ita and then set }
 \vskip -.15in
$$
 \BF_{km,kn}\ses \sum_{\la\part k}c_\la 
 \prod_{i=1}^{\l(\la)}Q_{ \la_im, \la_in}.
 \eqno 3.4
 $$ 
 
  \noindent
 The commutativity of the operators $Q_{u_1,v_1}$  and 
 $Q_{u_2,v_2}$  with $(u_1,v_1)$ and $(u_2,v_2)$ collinear 
 vectors assures that 3.4  well  defines the operator $\BF_{km,kn} $.
  In [3] a variety of Shuffle conjectures were formulated based on the above construction and various choices of the symmetric function $F[x]$. The simplest one corresponds to choosing $F=e_k$ (the $k^{th}$ elementary symmetric function). The corresponding operator which we denote $\Be_{km,kn}$ has a truly remarkable connection to the Theory of Parking Functions. 

\hsize 5.3 in 
 A single example will suffice  to get across this 
 connection. We have displayed here on the right a $12\times 20$ lattice rectangle with a path that proceeds  by  North and East unit steps, always remaining weakly above the diagonal $(0,0)\RA (12,20)$. A Parking Function in the $12\times 20$ lattice rectangle is obtained 
by labeling the cells immediately to the right of the north steps of such a path by the integers $1,2,\ldots ,20$ (referred to as  {\ita cars}) in a column increasing manner. This given, one of the conjectures formulated in [3] may expressed as the  identity

\vskip -.14in
$$
\Be_{km,kn}(-1)^{k(n+1)}\ses \sum_{PF\in {\cal PF}_{km,kn}}
t^{area(PF)} q^{dinv(PF)}F_{pides(PF)} 
\eqno 3.5
$$ 
 \vskip -2.4in
  \hsize 6.5 in
 \sapp
\hfill $
 {  \includegraphics[height=1.6 in]{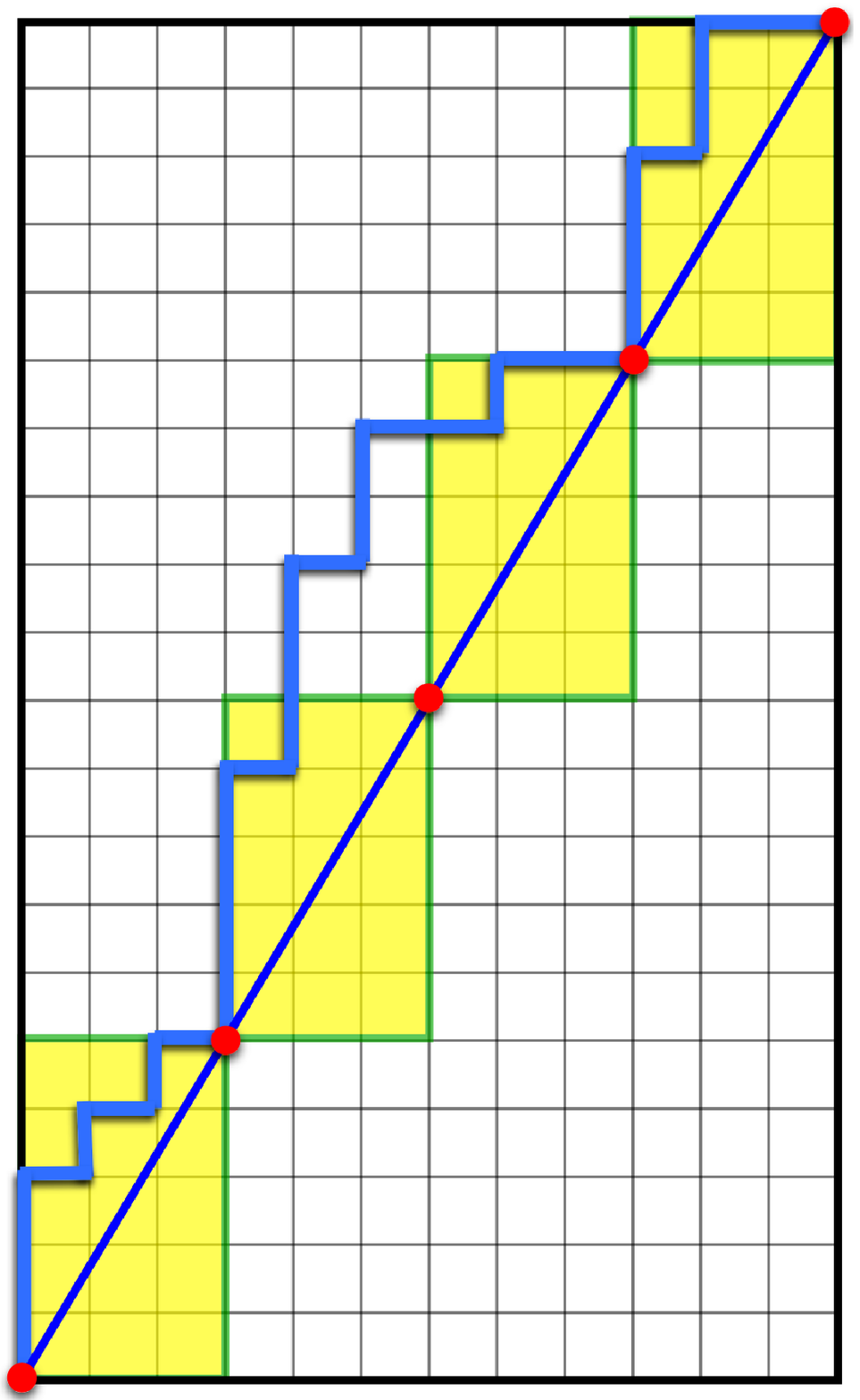} }
$
\vskip .13in

\noindent
where the sum is over all Parking Functions in the $km\times kn$ lattice rectangle, and the Parking Function statistics occurring in 3.5 are as 
will be defined shortly. Here the symbol ``$F_{pides(PF)} $'' stands for the
Gessel's [15] fundamental quasi-symmetric function indexed by the composition
``$pides(PF)$''.

Now the Extended Compositional Shuffle Conjecture  states a refined version of 3.5 where the symmetric function $F$ 
in the above construction of the operator $\BF_{km,kn}$, is chosen to assure that the sum in 3.5 is restricted  to be over an appropriately selected subfamily of Parking Functions in the   $km\times kn$ lattice rectangle.  These special choices of $F$ are obtained by means of the modified Hall-Littlewood operators $C_a$ whose action on a symmetric function $F[x]$ is defined by setting
\vskip -.08in
$$ 
 C_a P[X]= (-{\tttt{ 1\over q }})^{a-1} P\big[X-\tttt{1-1/q \over z}
\big]\sum_{m\ge 0}z^m h_m[X]\, \Big|_{z^a}
\eqno 3.6 
$$
\vskip -.06in

\noindent
Again, a single example will suffice to illustrate our choices.
 Suppose that we want to restrict the sum in 3.5 to be carried out only over the Parking Functions whose supporting Dyck path 
hits the diagonal precisely at the first and third and fourth possible places,
as the path depicted  in the above display. To achieve this we simply 
choose the symmetric function
 $
F=C_1C_2C_1\, 1
 $.
More generally, given a composition $p=(p_1,p_2,\ldots ,p_\l)$ of the integer $k$, let us denote by
$
\BC^{(p)}_{km,kn} 
$
the operator obtained by choosing $F=C_{p_1}C_{p_2}\cdots C_{p_\l}\, 1$ in the above construction. This given, the Extended Compositional Shuffle Conjecture states that
$$
\BC^{(p)}_{km,kn}(-1)^{k(n+1)}\ses \sum_{PF\in {\cal PF}_{km,kn}(p)}
t^{area(PF)} q^{area(PF)}F_{pides(PF)}
\eqno 3.7
$$ 
where the sum is over Parking Functions whose path hits the  diagonal precisely in $\l$ of the $k$ possible places,  as   prescribed  by the parts of the composition $p=(p_1,p_2,\ldots ,p_\l)$. The reason
3.7  refines 3.6 is due to the remarkable identity 
$$
e_k\ses \sum_{p\models k} C_{p_1}C_{p_1}\cdots C_{p_{\l(p)}}\, 1.
$$
Here the sum is over all compositions of $k$
\sas

Keeping all this in mind,  we are now in a position  to show that the identity in I.17 is one of the many consequences of the identity in 3.7.  More precisely we will show that  Theorem I.4 is a corollary of the following stronger result.
\sas

\noindent{\bol Theorem 3.1}

{\ita The validity of 3.7 for any coprime pair $(m,n)$, $k\ge1$ and any composition $p\models k$  implies the identity
$$
Q_{km,kn}(-1)^{kn} = \sum_{PF \in {\cal PF}_{km,kn}} t^{area(PF)-ret(PF)+1} \left[ ret(PF) \right]_t q^{dinv(PF)}
 F_{pides(PF)}
\eqno 3.8
$$
where all the Parking Function statistics are as  in the Extended Shuffle Conjectures. The ``$ret(PF)$'' statistic gives the smallest positive $i$ such that the supporting path of $PF$   goes through the point $(im,in)$.
}

\noindent{\bol Proof}

For brevity we will start with the following identity, valid for any integer $1\le   d\le k$ and  $p\models k-d$
$$
\sum_{p\models k-d; } \ssp C_d C _{p_1}C_{p_2}\cdots C_{p_{l(p)}}{\bf 1}
\ses (-{\TS{1\over q}})^{d-1}s_{d,1^{k-d}}- (-{\TS{1\over q}})^{d}\, 
     s_{1+d,1^{k-d-1}}  
\eqno 3.9
$$
(see [21] for a proof). Note next that it follows from 3.9  that
$$
 \sum_{d = a}^k \sum_{p\models k-d; } \ssp C_d    C_{p}\ 1
=
  \sum_{d = a}^k 
  (-{\TS{1\over q}})^{d-1}s_{d,1^{k-d}}
  -
 \sum_{d = a+1} ^k  (-{\TS{1\over q}})^{d-1}\, 
     s_{ d,1^{k-d }}  
= (-{\TS{1\over q}})^{a-1}s_{a,1^{k-a}}
 \eqno 3.10
$$
Our next step is to rewrite 3.1 in a more suitable form. Now a use of the Cauchy  identity gives
$$
\eqalign{
Q_{0,k} &=  
\tttt{(qt)^{1-k}\over qt-1}
\,  \ul h_{k}[X(1-qt)]
=
\tttt{(qt)^{1-k}\over qt-1}\sum_{\mu\part k}s_\mu[X]s_\mu[(1-qt]
\cr
&=
\tttt{(qt)^{1-k}\over qt-1}\sum_{r=0}^{k-1}
s_{k-r,1^r}[X](-qt)^r(1-qt)
=
-\tttt{(qt)^{1-k} }\sum_{r=0}^{k-1}
s_{k-r,1^r}[X](-qt)^r 
\cr
&=(-1)^k \sum_{a=1}^k s_{a,1^{k-a}}[X](-qt)^{1-a},
\cr}
\eqno 3.11
$$
where for the third equality  we use the following, easily verified,  special evaluation of a Schur function 
$$
s_\mu[1-m]\ses \cases{
(-m)^a(1-m) & if $\mu=k-r,1^r$ for $0\le r\le k-1$ \cr 
0 &  otherwise},
\eqno 3.12
$$ 
valid for any monomial $m$.

Hence by combining 3.10 and 3.11 we derive that
$$
\eqalign
{
(-1)^kQ_{0,k}  &=  \sum_{a=1}^k (-qt)^{1-a} s_{a,1^{k-a}} 
= \sum_{a=1}^k t^{1-a} \sum_{d=1 }^k \chi(d \geq a) \sum_{p\models k-d} C_d C_{p} 1 
\cr
&= \sum_{d=1}^k \sum_{p \models k-d} C_d C_{p} 1
 \sum_{a=1}^d  (1/t)^{a-1} 
= \sum_{d=1}^k \sum_{p \models k-d} [d]_{1/t} \, C_d C_{p} 1
\cr
}
\eqno 3.13
$$
This given, the particular case $F=Q_{0,k}$ of the above construction gives that for any coprime pair $(m,n)$ 
we have
$$
{\bf Q_{0,k}}_{km,kn}\ses Q_{km,kn}
\eqno 3.14
$$
Likewise by choosing $F=C_d C_{p} 1$ for $p\models k-d$ we obtain the operator 
 $
{\bf  C_{d,p} 1}_{km,kn}
$
which, by the Rational Shuffle Compositional conjecture, satisfies
$$
{\bf  C_{d,p} 1}_{km,kn}(-1)^{k(n+1)}
=\bu\bu
 \sum_{PF \in {\cal PF}_{km,kn}(d)}\bu\bu  t^{area(PF)}  q^{dinv(PF)}  F_{pides(PF)}
\eqno 3.15
$$
where the sum is over Parking Functions in the $km\times kn$ rectangle  whose supporting  path returns to  the diagonal for the 
first time in row $dn$. Thus combining 3.15 with 3.14 and 3.13 we obtain that
$$
Q_{km,kn}(-1)^{kn}
\ses \sum_{d=1}^k \sum_{p \models k-d} [d]_{1/t} \bu\bu\bu
 \sum_{PF \in {\cal PF}_{km,kn}(d)}\bu\bu  t^{area(PF)}  q^{dinv(PF)}  F_{pides(PF)}
$$
which is only another way of writing 3.8.
\sas
\hsize 5in
We can now finally give our 

\noindent{\bf  Proof of Theorem I.4}

We  will   show here the identity in I.17, written in the form
$$
\tttt{[k]_q\over [km]_q}     e_{kn}\big[ X[km]_q\big]
=\bu\bu\bu
\sum_{PF\in {\cal PF}_{km,kn}} \bu\bu\bu\bu
q^{coarea(PF)+dinv(PF)}
\big[ret(PF)\big]_q
F_{pides(PF)}[X]
\eqno 3.16
$$
\vskip -2.1in
  \hsize 6.5 in
 \sapp
\hfill $
 {  \includegraphics[height=1.8 in]{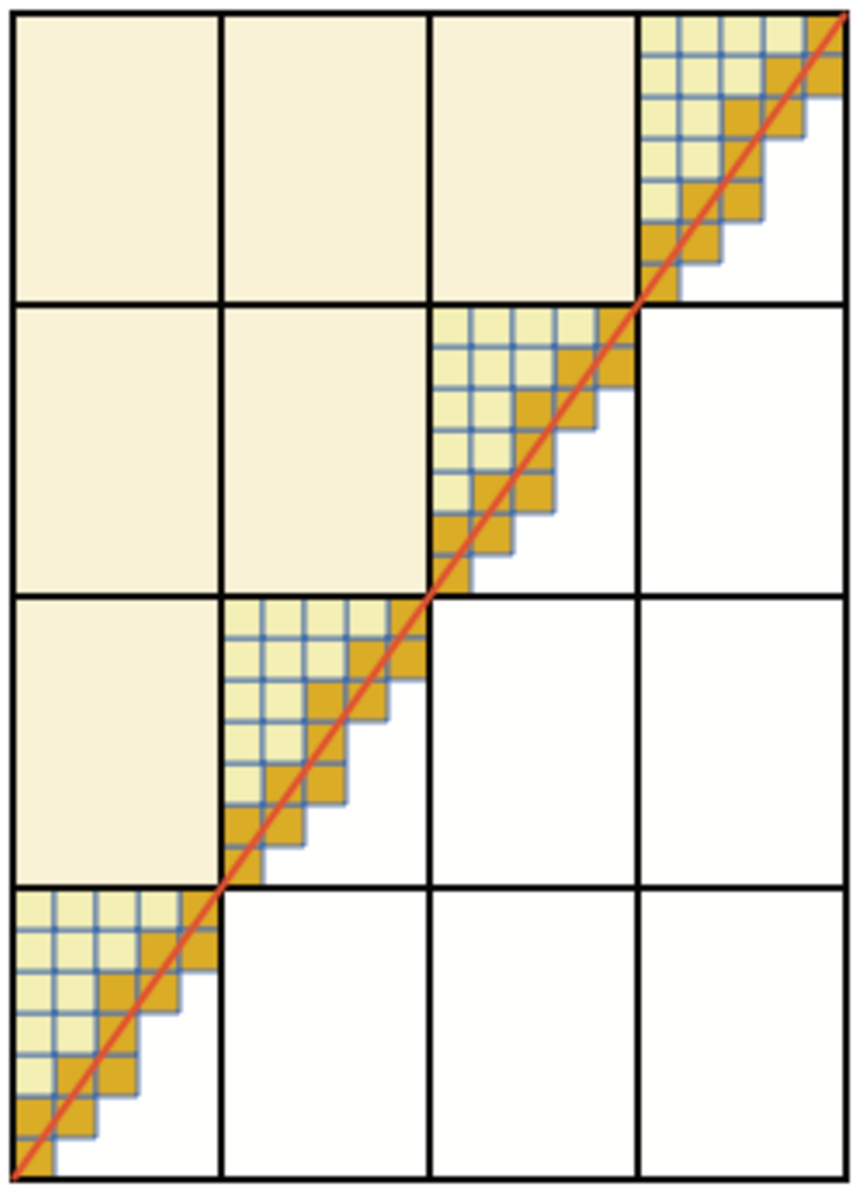} }
$

\vskip -.3in
\noindent
To  this end, notice that setting $t=1/q$ in 3.8 gives
$$
Q_{km,kn}(-1)^{kn}\Big|_{t=1/q} =  \bu\bu\bu\sum_{PF \in {\cal PF}_{km,kn}} \bu\bu\bu
q^{-area(PF)+ret(PF)-1 } \left[ ret(PF) \right]_q q^{-(ret(PF)-1)} q^{dinv(PF)}  F_{pides(PF)}
\eqno 3.17
$$

\noindent
But we can now use I.15 which can be  rewritten in the form 
$$
\tttt{[k]_q\over [km]_q}     e_{kn}\big[ X[km]_q\big]=
q^{(km-1)(kn-1)/2+( k-1)/2}  Q_{km,kn}(-1)^{kn} \Big|_{t=1/q}
\eqno 3.18
$$ 
This given, a comparison of the combination of 3.18 and 3.17 with 3.16 shows that we need only to show the equality
$$
coarea(PF)+area(PF)=(km-1)(kn-1)/2+( k-1)/2
\eqno 3.19
$$
This  can be easily justified by the following geometric argument. We have depicted in the display above the 
$km\times kn$ lattice rectangle for $k=4$ and $(m,n)=(5,7)$.
Now by the definition $area(PF)$ gives the number of lattice cells below the supporting path of $PF$ and weakly above 
the diagonal $(0,0)\RA (km,kn)$
and $coarea(PF)$ gives the number of lattice cells cells  above  the path. Thus to show 3.19 we need only verify that the right hand side gives he number of lattice cells  weakly above 
the diagonal $(0,0)\RA (km,kn)$. Now it is easy to see from the display that the number of lattice cells cut by the diagonal in
any one of the diagonal $4\times 7$ blocks is by
 one short of $4+7$ . This implies that the number of uncut lattice cells above the diagonal within each diagonal block 
 is none other than $(mn-m-n+1)/2=(m-1)(n-1)/2$. Moreover,
 the total  number of lattice cells within the upper non-diagonal blocks is $mn\times {k\choose 2}$. Thus 3.19 is none other than a consequence of the equality
$$
k\times (m-1)(n-1)/2\sps mn\times {k\choose 2}
\ses(km-1)(kn-1)/2+( k-1)/2
$$
This completes our proof of 3.16.
\sas

In the remainder of this section we will present the latest version of the Parking Function statistics that occur in the various formulations of the rational Compositional Shuffle Conjecture.

Let $(m,n)$ be coprime pair of positive integers and let $k \geq 1$. Recall that a $km,kn$-Dyck path is sequence of north and east steps in the $km \times kn$ lattice rectangle  which starts in the southwest corner, ends in the northeast corner, and stays weakly above the main diagonal $y = {n\over m} x$. A $km,kn$-parking function is a $km,kn$-Dyck path with labels $\{1,2,\dots,kn\}$, known as cars, adjacent to north steps and increasing within each column. For example, see the $6,9$-parking function below. Let ${\cal PF}_{km,kn}$ denote the set of all $km,kn$-parking functions.
\sas

\hsize 5 in
The original Rational Shuffle Conjecture of [19] states that for coprime $m,n$, we can express $Q_{m,n} (-1)^n$ as a weighted sum of $m,n$-parking functions. This enumeration involves the statistics $area(PF)$, $dinv(PF)$, and the word $\sigma(PF)$. The simplest of these is $area(PF)$, which is the number of full cells between the path and the main diagonal $y =  {n\over m} x$. In the adjacent example, the area is 5 and the corresponding cells are shaded.

The dinv and word statistics both make use of a rank function. Let $rank(x,y) = kmy - knx + \lfloor {x\over m} \rfloor$. This causes the points weakly above the diagonal to have distinct nonnegative ranks, with  points further from the main diagonal

\hsize 6.5in

 \vskip -1.9in
 \hfill $
 {  \includegraphics[width=1.3in]{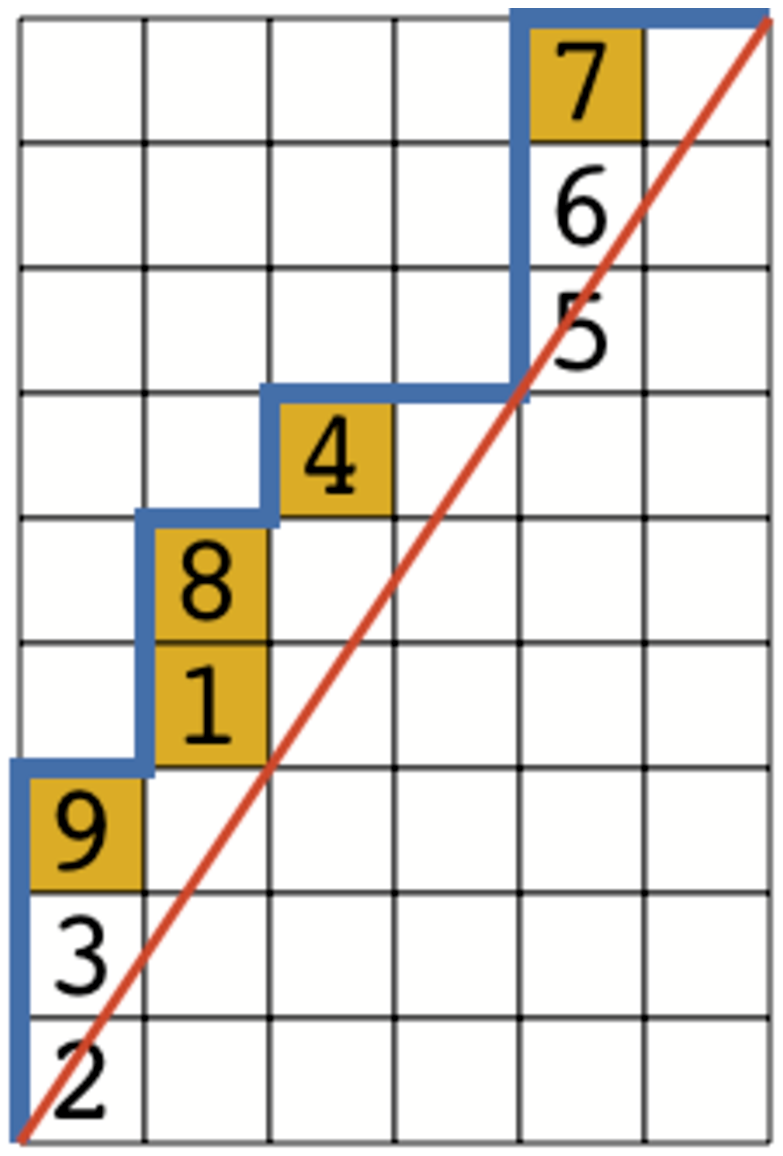} }
$
\vskip .01in

\vskip -.012in
\noindent
having higher rank.  This way, the rank function generalizes the notion of diagonals from classical Parking Function theory [20].
Let the rank of a car be the rank of the southwest corner of that car's cell. Then the word, $\sigma(PF)$, is just the permutation of $\{1,2,\dots,kn\}$ obtained by listing the cars from highest to lowest rank. In the example above, the ranks of cars $1$ through $9$ are $9$, $0$, $6$, $13$, $2$, $8$, $14$, $15$, and $12$, respectively. Hence the word of that parking function is $\sigma(PF) = 8 \, 7 \, 4 \, 9 \, 1 \, 6 \, 3 \, 5 \, 2$.

We set
$$
tdinv(PF) = \sum_{\hbox{cars }i<j} \chi(rank(i)<rank(j)<rank(i)+km).
\eqno 3.20
$$
Here $tdinv$ is short for "temporary dinv" because we will modify this statistic to obtain  $dinv(PF)$. In the example above, $tdinv(PF) = 9$ because the inequalities in 3.20 are satisfied for the pairs $(1,4)$, $(1,7)$, $(1,9)$, $(2,5)$, $(3,6)$, $(4,7)$, $(4,8)$, $(6,9)$ and $(7,8)$.

The original formulation of Hikita [24] as modified by Gorsky-Mazin [17], [18] expressed the dinv statistic as a combination of tdinv and two other statistics. However, Hicks and Leven [23] showed that this can be simplified as follows. Let $\lambda(PF)$ be partitions whose english Ferrers diagram is formed by the cells above $PF$. In the example above, $\lambda(PF) = (4,4,4,2,1,1)$. This given,   for an $m,n$-Parking  Function  
we set 
$$
dinv(PF) =
\cases{
 tdinv(PF) - \# \left\{ c \in \lambda(PF) : 
{arm(c)\over leg(c)} \leq {m\over n} <  {arm(c)+1\over leg(c)+1} \right\} & if $m<n$
\cr
  tdinv(PF) + \# \left\{ c \in \lambda(PF) :  {arm(c)\over leg(c)} >  {m\over n} \geq   {arm(c)+1\over leg(c)+1} \right\}
 & if $m>n$.
\cr
}
\eqno 3.21
$$ 
Here we must use the conventions  $ {0\over  0} = 0$ and $ {x\over 0} = \infty$ when  $x \neq 0$.
\supereject

\noindent
In the example above, we have $m<n$, thus
 $
 \ssp
dinv(PF) = tdinv(PF) - 4 = 5
\ssp  $,
since there are 3 cells with $arm=0$ and $leg=0$ and one cell with $arm=2$ and $leg=3$.

We now have all  the ingredients that occur in any of the
rational Shuffle  Conjectures including the Compositional ones in [3]. In particular,   3.5 conjectures the equality
$$
\Be_{km,kn}(-1)^{k(n+1)}\ses \sum_{PF\in {\cal PF}_{km,kn}}
t^{area(PF)} q^{dinv(PF)}F_{pides(PF)} 
\eqno 3.22
$$ 
where $pides(PF)$ is the composition that gives the descent set of the inverse of the permutation $\sig(PF)$ defined above,
and  $F_{pides(PF)}$ is the Gessel [15] fundamental quasi-symmetric function indexed by the composition $pides(PF)$.
Here we use the convention that  
 for a composition $p\models u$ if $S(p)$, is the subset of $\{1,2,\ldots , u-1\}$ that corresponds to $p$, then we set
$$
F_{p}(x_1,x_2,\ldots ,x_v)\ses
\sum_{\multi{1\le a_1\le a_2\le \cdots\le a_v\le v  \cr
i\in S(p)\RA a_i<a_{i+1} }} x_{a_1}x_{a_2}\cdots x_{a_v}
$$
The same conventions apply to the statistics occurring in 3.8,
namely  our conjecture that
$$
Q_{km,kn}(-1)^{kn} = \sum_{PF \in {\cal PF}_{km,kn}} t^{area(PF)-ret(PF)+1} \left[ ret(PF) \right]_t q^{dinv(PF)}
 F_{pides(PF)}
\eqno 3.23
$$
\sas

\noindent{\bol Remark 3.1}

We must mention that the conjectured equality in 3.16 has a 
 specialization that extends the equality in I.1 to the non coprime case. More precisely, we have
$$
\sum_{D\in {\cal D}_{km,kn}}[ret(PF)]_q  \ssp q^{coarea(D)+dinv(D)} \ses
{[k]_q\over [km]_q}\ssp \Big[ {kn+km-1\atop kn}\Big]_q
\eqno 3.24
$$

\hsize 5in
\noindent
In fact, scalar multiplication of both sides of 3.16 by   $e_{kn}[X] $
gives
$$
\eqalign{
\sum_{D\in {\cal D}_{km,kn}}[ret(PF)]_q  \ssp q^{coarea(D)+dinv(D)} &=
{[k]_q\over [km]_q}\ssp
\LL e_{kn}\big[ X[km]_q\big]\scs e_{kn}[X] \RR
\cr
&= 
{[k]_q\over [km]_q}\ssp \LL h_{kn}\big[ X[km]_q\big]\scs h_{kn}[X] \RR
\cr
}
$$
and 3.24 then follows from the identities
$$
\LL h_{kn}\big[ X[km]_q\big]\scs h_{kn}[X] \RR
= 
h_{kn}\big[[km]_q\big]=
\Big[ {kn+km-1\atop kn}\Big]_q
$$
\hsize 6.5in

 \vskip -2.1in
 \hfill $
 {  \includegraphics[width=1.6in]{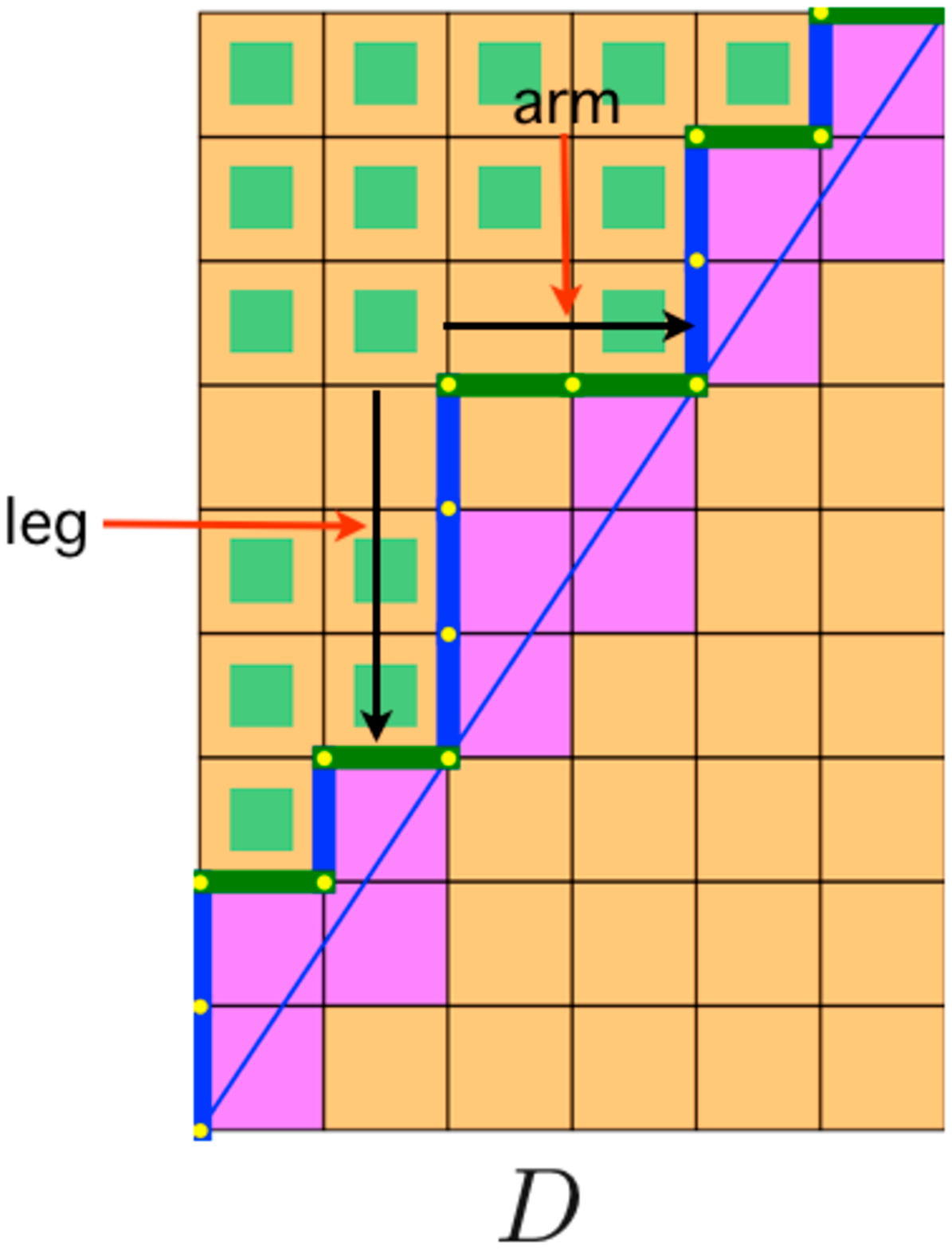} }
$

We should  mention that, for Dyck paths, the dinv statistic is obtained 
by counting the number of cells $c$ of the english partition above the path (see adjacent figure)  whose {\ita arm} and {\ita leg}  
satisfy the inequalities
$$
{arm(c)\over leg(c)+1}\le {m\over n}< {arm(c)+1\over leg(c)}
$$
In the example above we have placed a green square in each of the
cells that contribute to the dinv.
\sas

\noindent{\bol Remark 3.2}

This should complete  our presentation of the combinatorial side of the  rational Shuffle Conjectures except for two  
important observations. Firstly we should notice that 3.16
and I.17 differ in that I.16 has ``$s_{pides(PF)}$'' replacing
 ``$F_{pides(PF)}$'' in 3.16. We stated I.16 and various 
analogous  identities  in the introduction in this manner,
since this makes them  easier to verify on a computer.
In fact, the validity of this   replacement, 
is one of the  surprising  consequences  of  a  result  of  Egge, Loehr and Warrington [7]  concerning  Gessel  fundamental expansions of   symmetric functions. 

The second  observation results from a direct comparison of 3.22 an 3.23.  
Notice that we have  $1 \leq ret(PF) \leq k$ since the path must end at the point $(km,kn)$. Furthermore, $area(PF) \geq ret(PF)-1$. This is because each time $PF$ fails to touch the point $(im,in)$, a cell must fall between the path and the main diagonal. Therefore all the powers of $t$ appearing in 3.23 when 
$ret(PF)>1$ are non-negative. It follows from this that the difference of the right and sides of 3.23 and 3.22 can be shown to  be positive  linear combination of LLT polynomials 
which have in turn been shown to be Schur positive. 
So another strong evidence supporting the validity of these conjectures is that computer data confirms the Schur positivity 
of the difference of the left hand sides of 3.23 and 3.22.
\sa 

\noindent{\bol 4. The action of the operators $\bf D_{u,v}$
on the basis $\bf \big\{ s_\mu\big[{X\over 1-q} \big]\big\}_\mu$}
\sa
Our main goals  in this section are the proofs of Theorems I.5 and I.6.
To carry this out we need  auxiliary notation and some preliminary identities which may initially  appear only remotely connected 
with these goals.

Our basic tool here is a  new variant of the Macdonald  operator ``$D_n^1$'' of [25]. We will denote it ``$R_v$''.
Its  action on a symmetric polynomial $P[X_n]$ is obtained by setting
$$
R_vP[X_n]\ses {1 \over \DD[X_n]}\ssp
\sum_{i=1} ^n
T_{x_i}^qx_i^v\DD[X_n]P[X_n],
\eqno 4.1
$$
where ``$T_{x_i}^q$'' is the linear operator which carries out the substitution $x_i\RA qx_i$. Here $\DD[X_n]$ is the Vandermonde determinant$$
\DD[X_n]= 
\prod_{1\le  i<j\le n}(x_i-x_j) = \sum_{\sig\in S_n} \eee(\sig)\ssp \sig x^{\delta},
\eqno 4.2
$$
with $\eee(\sig)$ the sign of $\sig$ and   $\delta=(n-1,n-2,\ldots ,1,0)$. The actual value  of $n$
in both 4.1 and 4.2 
is immaterial provided that we choose it greater than $v$ plus the degree of $P$. The following identity, which in particular shows that $R_v$ preserves symmetry, will play  a crucial role. 
\sas

\noindent{\bol Proposition  4.1}

{\it For any integral vector $\mu=(\mu_1\ge \mu_2\ge \dots
\ge \mu_n\ge0)$ we have  
$$
R_vs_\mu[X_n]\ses  
 \sum_{i=1}^n q^{v+\mu_i+n-i}
 s_{\mu+ve_i}[X_n]
 \eqno 4.3
 $$
 with $e_i$ the $n$ dimensional coordinate vector  with 
$i^{th}$ component equal to $1$.}

\noindent{\bol Proof}

Since by definition
$$
s_\mu[X_n]\ses {1 \over \DD[X_n]} \sum_{\sig \in S_n}
\eee(\sig)\sig x^{ \mu+ \delta }
$$
we may write
$$
\eqalign{
R_vs_\mu[X_n]
&=
 {1 \over \DD[X_n]}
 \sum_{i=1} ^n
T_{x_i}^qx_i^v\sum_{\sig\in S_n} \eee(\sig)\sig x^{ \mu+\delta }
\cr
&=
 {1 \over \DD[X_n]}
 \sum_{i=1} ^n
 \sum_{\sig\in S_n} \eee(\sig)
 T_{x_{\sig_i}}^qx_{\sig_i}^v \sig x^{ \mu+\delta }
=
 {1 \over \DD[X_n]}
 \sum_{i=1} ^n q^{v+\mu_i+n-i}
 \sum_{\sig\in S_n} \eee(\sig)
   \sig x^{ \mu+ve_i+\delta }
\cr}
$$
 This proves 4.3.
\sas

The next identity shows that $R_v$ may be given an expression that is similar to the one obtained by Macdonald for his $D_n^1$ operator.
\sas

\noindent{\bol Proposition 4.2}

{\ita We have 
$$
R_v \ses q^v\sum_{i=1}^n \ssp A_i(x;q)\ssp x_i^vT_{x_i}^q
\eqno 4.4
$$ 
 with} 
$$
A_i(x;q)=  \prod_{j\in [1,n];j\neq i}  \bu {q\ssp x_i-x_j\over x_i-x_j}\ess .
\eqno 4.5
$$ 
{\bol Proof}

The  definition in 4.1 may also be rewritten as
$$
R_vP[X_n]\ses  
q^v\sum_{i=1} ^n
\Big({1 \over \DD[X_n]} T_{x_i}^q\DD[X_n]\Big) \ssp x_i^v T_{x_i}^qP[X_n]
$$
This shows 4.4 with
$$
A_i(x,q)\ses {1\over \DD[X_n]}
  T_{x_i}^q \DD[X_n] .
$$
However, we see that 
$$
\eqalign{
 {1\over \DD[X_n]}
  T_{x_i}^q \DD[X_n] &= 
    \prod_{ r,s\neq i}
  {x_r-x_s\over x_r-x_s}
    \prod_{i<s }
   {qx_i-x_s\over x_i-x_s}    
     \prod_ {r<i}
  {x_r-qx_i\over x_r-x_i}
  \ses
  \bu\bu\bu
  \prod_{j\in [1,n];j\neq i}  \bu {q\ssp x_i-x_j\over x_i-x_j},
    \cr
}
$$
as desired.
\sas

One of the difficulties in using Macdonald's original definition of the operators $D_n^k$ stems  from the fact that formulas expressing
a symmetric polynomial $P[X_n]$ in terms of the variables themselves are quite impractical
for significant values of $n$. For this reason, neither the definition in 4.1
nor its alternate form in 4.4 are much help in computing the action  the operators $D_n^k$ when it matters. However, for any of the operators $R_v$, we do have
a plethystic formula which is computationally as well as theoretically  very convenient. To state and prove this result  we need some auxiliary identities
\sas

\noindent{\bol Proposition 4.3}

{\ita For any formal series $F(x)=\sum_{k\geq 0}  c_k  x^k$  and for all integers $v\geq 0$ we have}
$$
F(1/z)\ssp \OM[(q-1)zX_n] \Big|_{z^v}=
{\chi(v=0)\over q^n}F(0)\sps
 {q-1\over q^n} 
\sum_{i=1}^n\ssp \Bigl(\prod_{\multi{j=1; j\neq i}}    {qx_i-x_j\over x_i-x_j}\Bigr)
  (qx_i)^v\ssp F(qx_i) 
  \eqno 4.6
$$ 
\noindent{\bol Proof}

Starting from the partial fraction expansion
$$
\OM[(q-1)zX_n]=\prod_{i=1}^n
 { 1-zx_i
\over 1-zqx_i
}=
 {1\over q^n}\sps {q-1\over q^n}\ssp 
\sum_{i=1}^n\ssp \Bigl(\prod_{\multi{j=1; j\neq i}}^n  {qx_i-x_j\over x_i-x_j}\Bigr)
{1\over 1-z qx_i}\ess ,
$$
we get 
$$
\eqalign{
F(1/z)\ssp \OM[(q-1)zX_n] \Big|_{z^v} 
&= 
{F(1/z)\over q^n}\Big|_{z^v}\sps {q-1\over q^n}\ssp 
\sum_{i=1}^n\ssp \Bigl(\prod_{\multi{j=1; j\neq i}}^n  {qx_i-x_j\over x_i-x_j}\Bigr)
{1\over 1-z qx_i}F(1/z)\Big|_{z^v}
\cr
&= 
{\chi(v=0)\over q^n}\, F(0)\sps {q-1\over q^n}\ssp 
\sum_{i=1}^n\ssp \Bigl(\prod_{\multi{j=1; j\neq i}}^n  {qx_i-x_j\over x_i-x_j}\Bigr)
\Big(\sum_{r\ge 0}(zqx_i)^r \Big)\Big(\sum_{s\ge 0}F_s/z^s \Big)\Big|_{z^v}
\cr
&= 
{\chi(v=0)\over q^n}\, F(0)\sps {q-1\over q^n}\ssp 
\sum_{i=1}^n\ssp \Bigl(\prod_{\multi{j=1; j\neq i}}^n  {qx_i-x_j\over x_i-x_j}\Bigr)
\Big(\sum_{\multi{r-s=v ; r,s\ge 0}}( qx_i)^r  F_s  \Big) 
\cr
&= 
{\chi(v=0)\over q^n}\, F(0)\sps {q-1\over q^n}\ssp 
\sum_{i=1}^n\ssp \Bigl(\prod_{\multi{j=1; j\neq i}}^n  {qx_i-x_j\over x_i-x_j}\Bigr)
\Big(\sum_{  s\ge 0} (qx_i)^{s+v}  F_s  \Big) 
\cr}
$$
This proves our Proposition.
\sas

As a corollary we obtain the following basic identity.
\sas

\vbox{
\noindent{\bol Theorem 4.1}

{\ita For any symmetric polynomial $P[X_n]$ we have
$$
R_v P[X_n]
\ses 
{\chi(v=0) \over 1-q}P[X_n] \sps {q^{n+v}\over q-1}\ssp 
P\bigl[ X_n-(1-q)/z \bigr] 
\ssp \OM\bigl[(1-1/q) zX_n\bigr]\ssp \big |_{z^v} .
\eqno  4.7
$$
}}

\noindent{\bol Proof}

Notice that in view of 4.4 we can rewrite 4.1 in the form
$$
\eqalign{
R_v P[X_n]
&\ses q^v\sum_{i=1}^n \ssp A_i(x;q)\ssp x_i^v
P\big[X_n-(1-q)x_i\big] 
\cr
&\ses 
q^v \sum_{i=1}^n\ssp \Bigl(\prod_{ j=1, j\neq i}^n  {qx_i-x_j\over x_i-x_j}\Bigr)\ssp x_i^v
P\bigl[ X_n-(1-q)x_i \bigr].
\cr
}
\eqno 4.8
$$ 
On the other hand Proposition 4.3 
 for  $F(z)\RA F(z/q)$ yields
$$
F(1/qz)\ssp \OM[(q-1)zX_n] \Big|_{z^v}=
{\chi(v=0)\over q^n}F(0)\sps
 {q-1\over q^n} 
\sum_{i=1}^n\ssp \Bigl(\prod_{\multi{j=1;j\neq i}}^n  {qx_i-x_j\over x_i-x_j}\Bigr)
  (qx_i)^v\ssp F(x_i) 
$$ 
Or better
$$
q^v\sum_{i=1}^n\ssp \Bigl(\prod_{ j=1, j\neq i }^n  {qx_i-x_j\over x_i-x_j}\Bigr)
\ssp  x_i ^v F(x_i)\ses {\chi(v=0) \over 1-q}F(0) \sps {q^n\over q-1}\ssp F(1/qz)\ssp \OM\bigl[(q-1)zX_n\bigr]\ssp \big |_{z^v} \ess .
$$
Using this with $F(z)=P\bigl[ X_n-(1-q)z \bigr]$ and using 4.8 gives
$$
\eqalign{
R_v P[X_n]
&\ses 
q^v\sum_{i=1}^n\ssp \Bigl(\prod_{ j=1, j\neq i }^n  {qx_i-x_j\over x_i-x_j}\Bigr)
\ssp  x_i ^v P\bigl[ X_n-(1-q)x_i \bigr]
\cr
&\ses 
{\chi(v=0) \over 1-q}P[X_n] \sps {q^n\over q-1}\ssp 
P\bigl[ X_n-(1-q)/qz \bigr] 
\ssp \OM\bigl[(1-1/q)qzX_n\bigr]\ssp \big |_{z^v} .
\cr
}
\eqno 4.9
$$
Notice next  that for any two formal power series $A(z),B(z)$ 
we have the identity
$$
\eqalign{
A[1 /qz]B[ qz ]\Big|_{z^v}
&=\sum_{r,s}A_rB_s(\tttt{1\over qz})^r(qz )^s\Big|_{z^v}
=\sum_{r,s }\tttt{ (zq)^{ s-r}}A_rB_{s}\Big|_{z^v}
\cr
&=
\sum_{s-r=v }  q ^v A_rB_{s} 
=q ^v\sum_{r }   A_rB_{r+v} =q ^vA[1 / z]B[  z ]\Big|_{z^v}
\cr
}
$$ 
and thus 4.9 becomes 
$$
R_v P[X_n]
\ses 
{\chi(v=0) \over 1-q}P[X_n] \sps {q^{n+v}\over q-1}\ssp 
P\bigl[ X_n-(1-q)/z \bigr] 
\ssp \OM\bigl[(1-1/q) zX_n\bigr]\ssp \big |_{z^v} .
$$
This proves 4.7.
\sas
 
 In particular, setting  $P[X_n]=s_\mu[X_n]$ and using 4.3 we obtain
 $$
\sum_{i=1}^n q^{v+\mu_i+n-i}
 s_{\mu+ve_i}[X_n]
\ses 
{\chi(v=0) \over 1-q}\ssp s_\mu[X_n] \sps {q^{n+v}\over q-1}\ssp 
s_\mu\bigl[ X_n-(1-q)/z \bigr] 
\ssp \OM\bigl[(1-1/q) zX_n\bigr]\ssp \big |_{z^v} .
$$
or better yet
$$
q^{n+v} \ssp 
s_\mu\bigl[ X_n-(1-q)/z \bigr] 
\ssp \OM\bigl[1-1/q) zX_n\bigr]\ssp \big |_{z^v} 
=
\chi(v=0)  s_\mu[X_n]\sps (q-1)\sum_{i=1}^n q^{v+\mu_i+n-i}
 s_{\mu+ve_i}[X_n]
\eqno 4.10
$$	
  
 At this point it is more convenient  to separate the cases
 $v>0$ and $v=0$. We will begin with the following 
 immediate corollary of Theorem 4.1.
  \sas
 
 \noindent{\bol Proposition  4.4}

{\ita For any $u,v>0$ and any partition $\mu$ we have  
$$
q^{uv}s_\mu[X - (1-q^u)/z ]  
\ssp \OM[(1-q^{-u}) zX] \Big|_{z^v}
\ses
 (q^u-1)\sum_{i=1}^{  |\mu|+v} q^{u(p(\mu)_i +v-i)}
s_{p(\mu)+ve_i }[X]
\eqno 4.11
$$
where $p(\mu)$ is the weak composition of length $ |\mu|+v $ obtained by adjoining zeros to the parts of $\mu$ and $e_i$ is the 
$i^{th }$ coordinate vector of length  $ |\mu|+v $. 
}

 \noindent{\bol Proof}

For $v>0$ 4.10 can be rewritten in the form
 $$ 
s_\mu\bigl[ X_n-(1-q)/z \bigr] 
\ssp \OM\bigl[1-1/q) zX_n\bigr]\ssp \big |_{z^v} 
=
 (q-1)\sum_{i=1}^n q^{ \mu_i+v -i}
 s_{\mu+ve_i}[X_n]
$$	
and the replacement $q\RA q^u$ gives
$$
q^{uv}s_\mu[X_n - (1-q^u)/z ]  
\ssp \OM[(1-q^{-u}) zX_n] \Big|_{z^v}
\ses
 (q^u-1)\sum_{i=1}^{ n} q^{u (\mu_i +v-i)}
s_{p(\mu)+ve_i }[X_n]
$$
This given, 4.11 follows  since the Schur functions involved in this expression stabilize after $n\ge  |\mu|+v$.
\sas

Keeping this in mind let us recall that our goal here is to 
 work out the action of the operator $D_{u,v}$ on the basis
  $\big\{s_\mu[\tttt{X\over 1-q}]\big\}_\mu$. The following identity
provides the link that ties this goal with the identities in 4.7 and 4.11
\sas

 \noindent{\bol Proposition  4.5}
 
 {\ita Suppose that for some $v\ge 0$ we have
 $$
s_\mu[X - (1-q^u)/z ]  
\ssp \OM[(1-q^{-u}) zX ] \Big|_{z^v}
\ses  G[X;q]
\eqno 4.12 
$$
 then }
 $$
 D_{u,v}s_\mu\big[\tttt{X\over 1-q}\big ]\ses q^v G\big[\tttt{X\over 1-q};q\big]
 \eqno 4.13
 $$
\noindent{\bol Proof} 
 
 Notice  first the following sequence of equalities
 $$
 \eqalign{
q^v G[X;q]
 &\ses
 q^v s_\mu\big[\tttt{X -(1-q^u)/ z  }\big]
\OM\big[ \tttt{ z}( 1- q^{-u} )X \big]\Big|_{z^v} 
\cr
&\ses s_\mu\big[\tttt{X -(1-q^u)/qz  }\big]
\OM\big[ \tttt{qz}( 1- q^{-u} )X \big]\Big|_{z^v}
\cr
&\ses s_\mu\big[\tttt{X + \tttt{ (q-1)\over 1-q}(1-q^u)/qz  }\big]
\OM\big[-\tttt{qz}\tttt{ 1- q^{-u}\over q- 1}X(1-q)\big]\Big|_{z^v}
\cr
&\ses s_\mu\big[\tttt{X + \tttt{ (1-1/q)\over 1-q}(1-q^u)/z  }\big]
\OM\big[-\tttt{z}\tttt{ 1- q^{-u}\over 1- q^{-1}}X(1-q)\big]\Big|_{z^v}
\cr
 }
 $$
 Next  the replacement $X\RA {X\over 1-q}$ gives
 $$
 \eqalign{
q^v G\big[\tttt{X\over 1-q};q\big]
 &\ses s_\mu\big[\tttt{{X \over 1-q}+ \tttt{ (1-1/q)\over 1-q}(1-q^u)/z  }\big]
\OM\big[-\tttt{z}\tttt{ 1- q^{-u}\over 1- q^{-1}}X\big]\Big|_{z^v}
\cr
 &\ses s_\mu\big[\tttt{{X+   (1-1/q) (1-q^u)/z \over 1-q}  }\big]
\OM\big[-\tttt{z}\tttt{ 1- q^{-u}\over 1- q^{-1}}X \big]\Big|_{z^v}
 \cr
 &\ses s_\mu\big[\tttt{{X+ (1-q)  (1-1/q) [u]_q/z \over 1-q}  }\big]
\OM\big[-\tttt{z}\tttt{ 1- q^{-u}\over 1- q^{-1}}X \big]\Big|_{z^v}
 \cr }
\eqno 4.14
 $$
 Now recalling that by definition we have  (for $t=1/q$)
 $$
D_{u,v}F[X]= F\big[X+ \tttt{M[u]_q/z}\big]
\OM\big[ -z[u]_tX\big]
\Big|_{z^v},
$$
we see that 4.14 proves 4.13.
 \sas
 
 We are thus able to obtain our
 \sas

\noindent{\bol Proof of Theorem I.6}

By combining Propositions 4.4 and 4.5 with
$$
G[X;q]\ses
q^{-uv}(q^u-1)\sum_{i=1}^{  |\mu|+v} q^{u(p(\mu)_i +v-i)}
s_{p(\mu)+ve_i }[X]
$$
we obtain
$$
 D_{u,v}s_\mu\big[\tttt{X\over 1-q}\big ]\ses
 (q^u-1)\sum_{i=1}^{  |\mu|+v} q^{up(\mu)_i+v -ui }
s_{p(\mu)+ve_i }\big[\tttt{X\over 1-q}\big ]
$$
as desired.
\sas

Our next task is to take care of the case $v=0$ of  4.7.
This may be rewritten as
$$ 
{q^{n}}\ssp 
P\bigl[ X_n-(1-q)/z \bigr] 
\ssp \OM\bigl[1-1/q) zX_n\bigr]\ssp \big |_{z^0}=
P[X_n]\sms (1-q)R_0 P[X_n]
$$ 
Choosing $P[X_n]=s_\mu[X_n]$ and using 4.3 for $v=0$ we get
$$ 
{q^{n}}\ssp 
s_\mu\bigl[ X_n-(1-q)/z \bigr] 
\ssp \OM\bigl[1-1/q) zX_n\bigr]\ssp \big |_{z^0}=
\Big(1\sms (1-q)
 \sum_{i=1}^n q^{\mu_i+n-i}
 \Big)
 s_{\mu}[X_n]
\eqno  4.15
$$ 

 Our next step is to transform 4.15 into a relation which contains no explicit dependence on $n$. 
To this end,  recalling that we write  a partition of $n$ in the form
  $\mu=(\mu_1\geq\mu_2\geq \cdots \geq \mu_n\geq 0)$ we set
$$
B_\mu(q,t)\ses \sum_{\mu_i>0} \ssp t^{i-1} (1+q+\cdots +q^{\mu_i-1})\ses 
\sum_{i=1}^n \ssp t^{i-1} \ess \tttt{1-q^{\mu_i}\over 1-q\ \ \ }
\eqno  4.16
$$
This polynomial, usually called the {\ita biexponent generator of }$\mu$, 
plays an essential role in the Theory of Macdonald polynomials, 
in our case this results in the following basic identity:
\sas

\noindent{\bol Proposition  4.6}

$$
s_\mu  \bigl[X_n-{\textstyle {1-q\over z}}\bigr]
\OM\bigl[(1-q^{-1})zX_n\bigr]\ssp \Big |_{z^0}\ses 
\Bigl(1\sms (1-q^{-1})(1-q)B_\la(q, q^{-1})\Bigr)\ssp s_\mu\bigl[X_n]
\eqno  4.17
$$ 
{\bol Proof}

Changing $t$ into $1/q$ in 4.16 and multiplying both sides 
 by $q^n$ we can rewrite it  in the form
$$
q^n\Big(1-(1- q^{-1})(1-q)B_\mu(q,t)\Big)\ses 
  1 
\sms (1-q)  \sum_{i=1}^n    q^{\mu_i+n-i}  
\eqno  4.18
$$
Using this in 4.15 gives
$$
 q^n s_\mu\bigl[ X_n-{\textstyle{1-q\over  z}} \bigr]
\OM\bigl[ (1-q^{-1}) zX_n\bigr]\ssp \big |_{z^0}=
q^n\Big(1-(1- q^{-1})(1-q)B_\mu(q,t)\Big)s_\mu[X_n]
$$
canceling the factor $q^n$ proves 4.17 as desired.
\sas

\noindent{\bol Remark 4.1}

By setting $v=0$ in 4.3 and 4.4 we derive that
$$
\sum_{i=1}^n \ssp A_i(x;q)\ssp  T_{x_i}^qs_\mu[X_n]\ses 
 \Big(\sum_{i=1}^n q^{ \mu_i+n-i}\Big).
 s_{\mu}[X_n]
$$
This identity was used by Macdonald in [25] to prove that his polynomial $P_\mu(X_n,q,t)$ reduces to $s_\mu[X_n]$
at $t=q$. Now the original definition of the modified 
Macdonald polynomial $\TH_\mu[X;q,t]$ was obtained by setting 
$$
\TH_\mu[X;q,t]\ses c_\nu(q,t)P_\mu[X/(1-1/t);q,1/t]
\eqno  4.19
$$
 where  $P_\mu[X;q,t]$ is none other than $P_\mu[X_n;q,t]$
(for $\mu\part n$) with $X_n$ replaced by the infinite alphabet
$X=x_1+x_2+\cdots$. and $c_\nu(q,t)$ is a polynomial in $q,t$
whose nature is immaterial here.
Thus it follows  from 4.19 that
$$
\TH_\mu[X;q,1/q]= c_\nu(q,q^{-1})P_\mu[X/(1-q);q,q]
=c_\nu(q,q^{-1})s_\mu \big[\tttt{X\over 1-q}\big]
 $$
Thus for all practical purposes, in the present context,
which arises from our setting $t=1/q$ in all our $Q_{u,v}$
and their identities, the basis $\{\TH_\mu[X;q,t]\}_\mu$ need 
only be replaced by the basis $\big\{s_\mu \big[\tttt{X\over 1-q}\big]\big\}_\mu$. In this vein we can easily obtain an alternate way
of interpreting the identity in 4.17 .
\sap

\noindent{\bol Theorem I.5}

{\ita For all integers $u\ge 1$ we have}
$$
D_{u,0}s_\mu\big[ \tttt{X\over 1-q} \big]\ses
\Bigl(1\sms (1-q^{-u})(1-q^u)B_\mu(q^u, q^{-u})\Bigr) s_\mu\big[ \tttt{X\over 1-q} \big]
\eqno  4.20
$$
\noindent{\bol  Proof}

Making the replacements $X_n\RA X$ and $q\RA q^u$ in 4.16 gives
$$
s_\mu  \bigl[X -{\textstyle {1-q^u\over z}}\bigr]
\OM\bigl[(1-q^{-u})zX \bigr]\ssp \Big |_{z^0}\ses 
\Bigl(1\sms (1-q^{-u})(1-q^u)B_\mu(q^u, q^{-u})\Bigr)\ssp s_\mu\bigl[X ]
$$ 
Next we do  $X \RA{X\over 1-q}$ and get
$$
s_\mu  \bigl[\tttt{X\over 1-q} -{\textstyle {1-q^u\over z}}\bigr]
\OM\bigl[  -\tttt{ 1-q^{-u} \over 1-1/q} (z/q)X\bigr]\ssp \Big |_{z^0}\ses 
\Bigl(1\sms (1-q^{-u})(1-q^u)B_\mu(q^u, q^{-u})\Bigr)\ssp s_\mu\bigl[\tttt{X\over 1-q} ]
$$ 
and this (with $t=1/q$) may be rewritten as 
$$
s_\mu  \bigl[
\tttt{
{X + {(1-1/q)(1-q) [m]_q/(z/q)}\over 1-q}
}\bigr]
\OM\bigl[-[m]_t(z/q)X \bigr]\ssp \Big |_{z^0}\ses 
\Bigl(1\sms (1-q^{-u})(1-q^u)B_\mu(q^u, q^{-u})\Bigr)\ssp s_\mu\bigl[\tttt{X\over 1-q} ]
$$ 
or better
$$
s_\mu  \bigl[
\tttt{
{X + {M [m]_q/z}\over 1-q}
}\bigr]
\OM\bigl[-[m]_tzX \bigr]\ssp \Big |_{z^0}\ses 
\Bigl(1\sms (1-q^{-u})(1-q^u)B_\mu(q^u, q^{-u})\Bigr)\ssp s_\mu\bigl[\tttt{X\over 1-q} ]
\eqno 4.21
$$ 
Recalling that by definition we have
$$
D_{u,v}F[X]= F\big[X+ \tttt{M[u]_q/z}\big]
\OM\bigl[-[u]_tzX \bigr]\ssp \Big |_{z^v}.
$$
We see that 4.20 is simply another way of writing  4.21.
This completes our proof.
\sas

\noindent{\bol Remark 4.2}

If we follow the sequence of steps that yielded the identity in 4.20 we will notice that this identity is but a direct consequence of the identity in 4.17
with the replacement $X_n\RA X$, that is  
$$
s_\mu  \bigl[X -{\textstyle {1-q\over z}}\bigr]
\OM\bigl[(1-q^{-1})zX \bigr]\ssp \Big |_{z^0}\ses 
\Bigl(1\sms (1-q^{-1})(1-q)B_\mu(q, q^{-1})\Bigr)\ssp s_\mu  [X ].
\eqno  4.22
$$ 
Thus in principle, a short cut in the proof of Theorem I.5, may  
appear to be the verification of 4.22.

 To better appreciate the power of the path we followed in the proof of 4.20, it will be instructive to see what kind of combinatorial identities  we are led to in trying  to carry this out.
Now working first  on the left hand side of 4.22 gives
$$
\eqalign{
LHS
& \ses
 s_\mu[X]+\sum_{k\ge 1}
\sum_{\nu =(k-a,1^{a})}s_{\mu/\nu'}[X] (-1/z)^k 
s_\nu[1-q]\OM\bigl[(1-q^{-1})zX \bigr]\ssp \Big |_{z^0}
\cr
& \ses
s_\mu[X]+(1-q)\sum_{k\ge 1}(-1)^k
\sum_{a=0}^{k-1}
s_{\mu/ (a+1,1^{k-a-1})}[X](-q)^ah_k[1-q^{-1})X]
\cr
}
\eqno 4.23
$$
and since
$$
h_k[1-q^{-1})X]\ses(1-q^{-1})\sum_{b=0}^{k-1}(-q)^bs_{k-b,1^b}[X]
$$
4.23 becomes
$$
LHS=s_\mu[X]+(1-q)(1-q^{-1})
\sum_{k\ge 1}(-1)^k
\sum_{a,b=0}^{k-1}(-q)^{q+b}s_{\mu/(a+1,1^{k-a-1})}[X]s_{k-b,1^b}[X]
\eqno 4.24
$$
Taking the scalar product of both sides of 4.22 by $s_\la[X]$ and using
4.24, routine manipulations reduce  4.22 to the  equivalent   identity
$$
\sum_{k\ge 1}(-1)^k
\sum_{a,b=0}^{k-1}(-q)^{a+b}
\LL s_{\mu/(a+1,1^{k-a-1})} \scs s_{\la/(k-b,1^b)} \RR
\ses -\chi(\la=\mu)B_\mu(q,q^{-1})
\eqno 4.25
$$
A standard result  on scalar products of skew Schur functions
(see [13]),  asserts  that the scalar product summand  is none other than  
 the number of permutations that fit the shape 
$\mu/(a+1,1^{k-a-1})$ whose inverse fits the shape $\la/(k-b,1^b)$.
This given, the only conclusion we can draw from this calculation   is that 4.25 is one truly remarkable combinatorial consequence of  Theorem I.5. 
\page

\noindent{\bol 5. The original proof of Theorem I.1 by the partial fraction method}

\sas

In this section we explain how we discovered Theorem I.1. Indeed Theorem I.1 involves  different series expansions of a single rational function,
which is best understood by using the partial  fraction method of the fourth named author.
To this end we need to work in the field $K=\BQ((z_N))((z_{N-1}))\cdots ((z_1))$ of iterated Laurent series to obtain series expansion of rational functions.
The readers are referred to [31] for the original development of the field of iterated Laurent series. Here we only recall that $K$ defines a total group order on its monomials given by
$$
  z_1^{a_1}\cdots z_N^{a_N}  
                              \cases{ 
                                <_K 1, &  if   $a_1=\cdots =a_{i-1}=0 \; \& \; a_{i}>0$;  \cr 
                                =1, &  if $ a_1=\cdots =a_N=0$;
\cr 
                                >_K 1, &    if  $a_1=\cdots =a_{i-1}=0 \; \& \; a_{i}<0$.
\cr
}
$$
We shall simply write this order by $z_1<z_2<\cdots <z_N<1$. The series expansion of $(1-w)^{-1}$ for a monomial $w\ne 1$ (called small or large) is thus given by
$$ 
 {1\over1-w} = \cases{
                    { \sum_{n\ge 0}}  w^n, &   if   $w<_K 1$; 
                     \cr\cr 
                     {1\over -w(1-1/w)}=-\sum_{n\ge 0} w^{-n-1}, &  if  $ w>_K 1$.
}
 $$
To start with, let us recall
that for any Laurent polynomials $L(z_1,z_2)\in \Q[z_1^{\pm 1},z_2^{\pm 1}]$,
we always have that
$$ 
 L(z_1,z_2) \Big|_{z_1^0z_2^0}=L(z_2,z_1) \Big|_{z_1^0z_2^0}    \qquad \hbox{holds in }  \Q[z_1^{\pm 1},z_2^{\pm 1}]. 
 $$
In particular, antisymmetric Laurent polynomials have constant term $0$.
Such properties no longer hold for rational functions in $K_1=\Q((z_1))((z_2))$. For example,
$f(z_1,z_2)=
 {z_1+z_2\over z_1-z_2}$ is clearly antisymmetric, but in  $K_1$ we have
$$ 
 {z_1+z_2\over z_1-z_2}\Big|_{z_1^0z_2^0} =  {1+z_2/z_1\over 1-z_2/z_1} \Big|_{z_1^0z_2^0}= (1+z_2/z_1) \sum_{n\ge 0} (z_2/z_1)^n \Big|_{z_1^0z_2^0} =1.
$$
Indeed, the exchanging of the two variables $z_1,z_2$ transforms the constant term in $K_1$ to a constant term in $K_2=\Q((z_2))((z_1))$. We shall have
$$
 f(z_1,z_2) \Big|^{K_1}_{z_1^0z_2^0} = f(z_2,z_1) \Big|^{K_2}_{z_1^0z_2^0}
 $$
where the left hand side is a constant term in $K_1$, but the right hand side is a constant term in $K_2$. This type of exchanging of variables would be sufficient for us to prove Theorem I.1. See [30]
 for general formulation on the change of variables in a field of Malcev-Neumann series.
\sas

\noindent{\bol Proof of Theorem I.1}

With $t=1/q$, we have
$$
\eqalign{
  D_{c,d} D_{a,b} F[X] &= D_{c,d} F\big[X+ [a]_q \tttt   {M\over z_1}\big]
  \Omega\big[-z_1 X[a]_t\big] z_1^{-b} \Big|_{z_1^0} 
  \cr
&=
F\big   [X+ [c]_q   \tttt {M\over z_2}+ [a]_q   \tttt{M\over z_1}\big   ] \Omega\big   [-z_1 [a]_t\big   (X+[c]_q  \tttt {M\over z_2}\big   )\big   ] \Omega\big   [-z_2 X [c]_t\big   ] z_1^{-b} z_2^{-d} \Big|_{z_1^0 z_2^0}  \cr
&=
F\big   [X+ [c]_q    \tttt{M\over z_2} + [a]_q    \tttt{M\over z_1}\big   ] \Omega\big   [-X(z_1 [a]_t+z_2 [c]_t)\big   ] \Omega\big   [- M[a]_t [c]_q    \tttt{z_1\over z_2} \big   ]  z_1^{-b} z_2^{-d} \Big|_{z_1^0 z_2^0} 
\cr
&=
F\big   [X+ [c]_q   \tttt {M\over z_2} + [a]_q   \tttt {M\over z_1}\big   ] \Omega\big   [-X(z_1 [a]_t+z_2 [c]_t)\big] \displaystyle   {(1-   {z_1\over z_2})(1- q^ct^a   {z_1\over z_2})\over (1-q^c    \tttt{z_1\over z_2})(1-t^a   \tttt {z_1\over z_2})}  z_1^{-b} z_2^{-d} \Big|_{z_1^0 z_2^0}
 \cr
} 
$$
where in the last step we have used the fact that $M[a]_t[c]_q=(1-t^a)(1-q^c)=1-t^a-q^c+q^ct^a$. This constant term has to be understood as in a field of iterated Laurent series
where $q^c z_1/z_2$ and $t^a z_1/z_2$ are small. In the general $q,t$ case,
we can set $q<t<z_1<z_2<1$. But here we can not set $q$ to be small, since that will force $t=1/q$ to be large. We choose to work in the field of iterated Laurent series $K_1$ defined by the order $z_1<z_2<q<1$ (one can take $q$ as a constant). More precisely, we can set, for example, $K_1=\Q(q)((z_2))((z_1))[[x_1,x_2,\dots]]$.

Let us write
$$
   D_{c,d} D_{a,b} F[X] = F\left[X+ [c]_q \   \tttt  {M\over z_2} + [a]_q \    \tttt   {M
\over z_1}\right] \Omega\left[-X(z_1 [a]_t+z_2 [c]_t)\right] G(z_1,z_2;c,d,a,b)     \Big|^{K_1}_{z_1^0 z_2^0},
\eqno 5.1
$$
where
$$ G(z_1,z_2;c,d,a,b)= \     {(1-\   \tttt    {z_1\over z_2})(1- q^ct^a\   \tttt    {z_1
\over z_2})\over (1-q^c \   \tttt    {z_1\over z_2})(1-t^a  {z_1\over z_2}) }  \ssp z_1^{-b} z_2^{-d}.
$$
Now switching $(c,d)$ and $(a,b)$ gives
 $$ 
 D_{a,b} D_{c,d} F[X] = F\left[X+ [a]_q \    \tttt   {M\over z_2} + [c]_q 
 \  \tttt     {M\over z_1}\right] \Omega\left[-X(z_1 [c]_t+z_2 [a]_t)\right] G(z_1,z_2;a,b,c,d)     \Big|^{K_1}_{z_1^0 z_2^0}.
 $$
Observe that when $t=1/q$, we have
$$
G(z_2,z_1;a,b,c,d)= \     {(1-\    \tttt   {z_2\over z_1})(1- t^cq^a 
\  \tttt     {z_2\over z_1})\over (1-t^c \   \tttt    {z_2\over  z_1})(1-q^a \     \tttt  {z_2\over  z_1})}=\     {(1-\    \tttt  {z_1\over  z_2})(1- q^ct^a\   \tttt    {z_1\over z_2})\over  (1-q^c \   \tttt    {z_1\over  z_2})(1-t^a \  \tttt     {z_1\over  z_2})}=G(z_1,z_2;c,d,a,b).
$$

By exchanging $z_1$ and $z_2$, we obtain that $D_{a,b} D_{c,d}  F[X]$ is the same constant term as in 5.1, but working in the field of iterated Laurent series $K_2$
defined by the order $z_2<z_1<q<1$. So we are indeed computing the difference of the constant terms of a single ``rational function" in two different working fields.

By partial fraction decomposition in $z_1$, applied to the coefficients in the $x$'s and then sum, we have
$$
\eqalign{
F\left[X+ [c]_q  \tttt{M\over z_2} + [a]_q  \tttt{M\over z_1}\right]  &\Omega\left[-X(z_1 [a]_t+ z_2 [c]_t)\right] G(z_1,z_2;a,b,c,d) =
\cr
&\bigsp \bigsp = p^{\ge 0}  (z_1)+ p^{<0} (z_1)+ {A_1\over 1-q^c   {z_1\over z_2}} + {A_2\over 1-t^a  {z_1\over z_2}},
\cr
}
$$
where when we restrict to each coefficient of the $x$'s, $p^{\ge 0}(z_1)$ is a polynomial in $z_1$, $p^{< 0}(z_1)$ is a Laurent polynomial that only contains negative powers in $z_1$, and
$A_1$ and $A_2$ are free of $z_1$ given by

$$
\eqalign{
  A_1 &
=F\left[X+ [c]_q  \tttt{M\over z_2} + q^c [a]_q \     \tttt  {M\over z_2}\right] \Omega\left[-X(z_2 t^{c} [a]_t+z_2 [c]_t)\right] \       {(1-t^c)(1-t^a)\over (1-t^a t^c)}  z_2^{-b}t^{-bc} z_2^{-d}, 
\cr
&=-q^{bc}\       {(1-q^c)(1-q^a)\over (1-q^{a+c})}   F\left[X+ [a+c]_q
\      \tttt {M\over  z_2} \right] \Omega[-Xz_2 [a+c]_t ]   z_2^{-b-d}   ;
\cr
}
\eqno 5.2
$$ 
$$
\eqalign{
A_2=&F\left[X+ [c]_q \   \tttt    {M\over z_2}+ t^a[a]_q \    \tttt   {M\over z_2}\right] \Omega[-X(z_2 q^a [a]_t+z_2 [c]_t)] \       {(1-q^a)(1- q^c)\over  (1-q^a q^c)}  z_2^{-b}q^{-ba} z_2^{-d} 
\cr
=&\       {(1-q^a)(1- q^c)\over (1-q^a q^c)} F\left[X+ [a+c]_q  \      \tttt {M\over z_2 q^a}\right] \Omega[-Xz_2 q^a[a+c]_t]   z_2^{-d-b} q^{-ba}.
\cr
}
$$ 
Actually we need the following formula (obtained by the change of variables $z_2\RA  z_2 t^a$).
$$
A_2\Big|_{z_2^0}= q^{da}\       {(1-q^a)(1- q^c)\over (1-q^{a+c})} F\left[X+ [a+c]_q  \       \tttt{M\over z_2}\right] \Omega[-Xz_2 [a+c]_t]   z_2^{-d-b} \Big|_{z_2^0}.  
\eqno 5.3
$$

Now we take the constant term in $z_1$ first, working in $K_1$ and $K_2$ separately. Since the two monomials $q^c z_1/z_2$ and $t^a z_1/z_2$ are both small in $K_1$ but large in $K_2$, we have
$$
\eqalign{
 D_{c,d} D_{a,b} F[X] 
 &= p^{\ge 0}(0) \Big|_{z_2^0} + A_1 \Big|_{z_2^0} + A_2 \Big|_{z_2^0}. \cr
D_{a,b} D_{c,d} F[X] &= p^{\ge 0}(0) \Big|_{z_2^0}.
\cr
}
$$ 
It follows that
$$ 
    {1\over M} [D_{c,d}, D_{a,b}] F[X] =     {1\over M} \Big(A_1 \Big|_{z_2^0} + A_2 \Big|_{z_2^0}\Big). 
$$
Applying formulas  5.2 and 5.3 gives the desired result.
\sas

The idea of the proof of Theorem I.1 can be generalized as follows, which is easy to prove but turns out to be very useful.
\sas

\noindent{\bol Proposition 5.1}

{\ita Let $K_1$ and $K_2$ be two different field of iterated Laurent
series. Suppose that we have the following partial fraction
expansion:
$$
F(z)=p_0(z)+{p_{-1}(z))\over  z^m}+ {p_1(z)\over (1-u_1z)^{k_1}}+\cdots +  {p_N(z))\over (1-u_Nz)^{k_N}}.
$$
Then we have 
$$ 
F(z)\Big|_{z^0}^{K_1}   - F(z)\Big|_{z^0}^{K_2} = \sum_{u_iz <_{K_1} 1\; \&\; u_iz >_{K_2} 1 } p_i(0) - \sum_{u_jz>_{K_1} 1\; \& \; u_jz<_{K_2} 1 }p_j(0). 
$$
In words, the difference of the two constant terms only came from
those denominators that are contributing (i.e., with $u_iz<1$) in one field but dually
contributing (i.e., with $u_iz>1$) in the other field.}

\noindent{\bol Proof}

With the given partial fraction decomposition, taking constant term in $z$ under $K_1$ gives
$$
F(z) \Big|_{z^0}^{K_1} = p_0(0) +\sum_{u_iz^{k_i}<_{K_1} 1} p_i(0).
$$
A similar result holds for $K_2$. Subtracting gives the desired formula.
\sas

Remark: The proposition  applies if $F(z)$ includes something like
$\Omega[X+M/z]$ or $\Omega[-zX]$ as factors, in which case $m$ or
$p_0(z)$ does not exist. But $F(z)$ can be first expanded as a power
series in the $x$'s, and then apply the proposition to the
coefficients in the $x$'s.

Remark: If we take $K_1={\BC}((z))$ and $K_2= {\BC}((z^{-1}))$, then the proposition gives
$$ F(z)\Big|_{z^0}^{K_1}   - F(z)\Big|_{z^0}^{K_2} = \sum_{i} p_i(0) . $$
This can be shown to be equivalent to the well-known fact that for any given rational function,
its residues at all points (including $\infty$)
sum to $0$.
\page

 \centerline{\bol Bibliography}
\sas
\item{[1]} F. Bergeron, A. M. Garsia, M. Haiman and G. Tesler, {\ita Identities and positivity conjectures for some remarkable operators in the theory of symmetric functions}, Methods Appl. Anal. 6 (1999), 363Ð420.
\sas

\item{[2]} F. Bergeron and A. M. Garsia,  {\ita Science Fiction and Macdonald's Polynomials}.
 CRM Proceedings  $\&$ Lecture Notes, American Mathematical
  Society , 22:  1--52, 1999.
\sas

\item{[3]} F. Bergeron, A.  Garsia, E. Leven, G. Xin,
{\ita Compositional (km,kn)-Shuffle Conjectures},
math arXiv:1404.4616.
\sas

\item{[4]}
 F. Bergeron, A. M. Garsia,  E. Leven and G. Xin,   
{\ita Some remarkable new Plethystic operators in the Theory of Macdonald Polynomials}, math arXiv: 1405.0316.
\sas

\item{[5]} I. Burban and O. Schiffmann,
{\ita  On the Hall algebra of an elliptic curve, I,} Preprint (2005), arXiv:math.AG/0505148.
\sas

\item{[6]} C. Dunkl, {\ita Intertwining operators and polynomials associated with the
symmetric group}, Monatsh. Math. 126 (1998), no. 3, 181Ð209.
\sas

\item{[7]} 
E. Egge, N. Loehr and G. Warrington,
{\ita From quasisymmetric expansions to Schur expansions via a modified inverse Kostka matrix},
European Journal of Combinatorics, {\BV}. 1. Dec 2010 pp. 2014-2027
\sas

\item{[8]} 
A. M. Garsia and J. Haglund, {\ita Proof of the $q,t$-Catalan Positivity Conjecture}, Discrete Mathematics, {\BV} 256, Issue 3, 28 October 2002, pp. 677-717
\sas

\item{[9]} 
A.M. Garsia, J. Haglund, G. Xin, M. Zabrocki
{\ita Some new applications of the Stanley-Macdonald Pieri Rules},  arXiv:1407.7916, Comments:  for Stanley@70. Massachusetts institute of technology, June 23-27, 2014
\sas

\item{[10]}
A.  Garsia and M. Haiman, {\ita A remarkable q, t-Catalan sequence and q-Lagrange inversion}, J. Algebraic
Combin. 5 (1996), no. 3, 191-244.
\sas

\item{[11]}
A. M. Garsia and  M. Haiman,
{\ita Some Natural Bigraded $S_n$-Modules and q,t-Kostka Coefficients. (1996) (electronic)
The Foata Festschrift}, http://www.combinatorics.org/Volume 3/volume 3 2.html\#R24
\sas

\item{[12]} 
A.~M. Garsia, M.~Haiman, and G.~Tesler.
 {\ita Explicit plethystic formulas for Macdonald $(q,t)$-Kostka
  coefficients}.
 S\'eminaire Lotharingien de Combinatoire [electronic only],
  42:\penalty0 B42m, 1999.
\sas

\item{[13]}A. Garsia and J. Remmel,  {\ita Shuffles of permutations and the Kronecker product}, Graphs and Combinatorics, V {\bf 1}, issue {\bf 1},
(1985) pp. 217-263.
\sas

\item{[14]} 
A. M. Garsia, G.~Xin and M. Zabrocki, 
 {\ita HallÐLittlewood Operators in the Theory of Parking Functions and Diagonal Harmonics
}, Int Math Res Notices (2011)
doi: 10.1093/imrn/rnr060
\sas

\item{[15]} 
I. Gessel.
 {\ita Multipartite P-partitions and inner products of skew Schur
  functions}.
 Contemp. Math, 34:\penalty0 289--301, 1984.
\sas

\item{[16]} E. Gorsky, {\ita Arc spaces and DAHSA representations}, math arXiv:1110.1674 (2011).
\sas

\item{[17]} E. Gorsky and M. Mazin, {\ita Compactified Jacobians and q,t-Catalan Numbers I},  
Journal of Combinatorial Theory, Series A, 120 (2013) 49-63.
\sas

 \item{[18]}, E. Gorsky and M. Mazin, {\ita Compactified Jacobians and q,t-Catalan Numbers II},  
Journal of Alg, Comb.l  (2014) {\bf 39}: 153-186.
\sas

\item{[19]} 
E. Gorsky and A.~Negut.
 {\ita Refined knot invariants and Hilbert schemes}.
 math arXiv:1304.3328, 2013.
\sas

\item{[20]} 
J. Haglund, M.~Haiman, N.~Loehr, J.~B. Remmel, and A.~Ulyanov.
 {\ita A combinatorial formula for the character of the diagonal
  coinvariants}.
 Duke J. Math., 126:\penalty0 195-232, 2005.
\sas

\item{[21]} 
J. Haglund, J. Morse, and M. Zabrocki.
 {\ita A compositional refinement of the shuffle conjecture specifying touch
  points of the Dyck path.}
 Canadian J. Math, 64:   822-844, 2012.
\sas

\item{[22]} M. Haiman, {\ita Conjectures on the Quotient Ring by Diagonal Invariants}, Journal of Alg. Comb. {\bf 3} (1994),17-76.
\sas

\item{[23]} A. Hicks and E. Leven. {\ita A simpler formula for the number of diagonal inversions of an (m,n)-Parking Function and a returning Fermionic formula },  math arXiv:1404.4889  (2014): (to appear in Discrete Math.) 
\sas

\item{[24]} 
T.~Hikita.
 {\ita Affine springer fibers of type A and combinatorics of diagonal
  coinvariants.}
 math arXiv:1203.5878, 2012.
\sas

\item{[25]}
I. G. Macdonald, {\ita Symmetric functions and Hall polynomials}, Oxford Mathematical Monographs, Oxford Science Publications, Oxford University Press, New York, 1995.
\sas

\item{[26]} 
O.~Schiffmann and E.~Vasserot.
 {\ita The elliptical Hall algebra, Cherednik Hecke algebras and Macdonald
  polynomials}.
 Compos. Math., 147.1:\penalty0 188--234, 2011.
\sas

\item{[27]} 
O.~Schiffmann and E.~Vasserot.
 {\ita The elliptical Hall algebra and the equivariant K-theory of the
  Hilbert scheme of ${A}^2$}.
 Duke J. Math., 162.2: 279--366, 2013.
\sas

\item{[28]} 
O. Schiffmann, {\ita On the Hall algebra of an elliptic curve, II}, Preprint (2005),
arXiv:math/0508553.
\sas

\item{ [29]} S. Stevan. Chern-Simons invariants of torus links. Ann. Henri Poincar«e 11 (2010), no. 7, 1201-1224.
\sas

  \item{[30]}
G. Xin,  {\ita A residue theorem for Malcev-Neumann series}, Adv. Appl. Math. (2005) 35, 271--293.
\sas

\item{[31]}
 G. Xin,  {\ita A fast algorithm for MacMahon's
partition analysis}, Electron. J. Combin., (2004), 11, R53.

\end

\end